\documentclass[sn-aps]{sn-jnl}

\usepackage{amssymb}
\usepackage{amsthm}
\usepackage{lineno}
\usepackage{wrapfig}
\usepackage{soul}
\usepackage[utf8]{inputenc}

\usepackage{amsmath}
\usepackage{url}
\usepackage{units}
\usepackage{algpseudocode}  
\usepackage{algorithm}      
\algnewcommand\algorithmicto{\textbf{to}}
\algnewcommand\algorithmicin{\textbf{in}}
\algnewcommand\algorithmicforeach{\textbf{for each}}
\algrenewtext{For}[3]{\algorithmicfor\ #1 $\gets$ #2\ \algorithmicto\ #3\ \algorithmicdo}
\algdef{S}[FOR]{ForEach}[2]{\algorithmicforeach\ #1\ \algorithmicin\ #2\ \algorithmicdo}

\renewcommand{\p}{\emph{p}}
\renewcommand{\r}{\emph{r}}
\newcommand{\h}{\emph{h}}
\newcommand{\hp}{\emph{hp}}
\newcommand{\hr}{\emph{hr}}
\renewcommand{\R}{\mathbb{R}}
\newcommand{\N}{\mathbb{N}}
\renewcommand{\L}{\mathcal{L}}
\renewcommand{\b}{\boldsymbol}
\newcommand{\lap}{\nabla^2}

\renewcommand{\x}{\b{x}}
\renewcommand{\w}{\b{w}}

\jyear{2023}%

\theoremstyle{thmstyleone}%
\theoremstyle{thmstyletwo}%

\theoremstyle{thmstylethree}%

\begin{document}

\title[Mesh-free \hp-adaptive solution procedure]{Strong form mesh-free \hp-adaptive solution of linear elasticity problem}


\author*[1,2]{\fnm{Mitja} \sur{Jančič}}\email{mitja.jancic@ijs.si}

\author[1]{\fnm{Gregor} \sur{Kosec}}\email{gregor.kosec@ijs.si}

\affil[1]{\orgdiv{Parallel and distributed systems laboratory}, \orgname{Institute Jožef Stefan}, \orgaddress{\street{Jamova Cesta 39}, \city{Ljubljana}, \postcode{1000}, \country{Slovenia}}}
\affil[2]{\orgname{Jožef Stefan International Postgraduate School}, \orgaddress{\street{Jamova Cesta 39}, \city{Ljubljana}, \postcode{1000}, \country{Slovenia}}}


\abstract{
    We present an algorithm for \hp-adaptive collocation-based mesh-free numerical analysis of partial differential equations. Our solution procedure follows a well-established iterative solve--estimate--mark--refine paradigm. The solve phase relies on the Radial Basis Function-generated Finite Differences (RBF-FD) using point clouds generated by advancing front node positioning algorithm that supports variable node density. In the estimate phase, we introduce an Implicit-Explicit (IMEX) error indicator, which assumes that the error relates to the difference between the implicitly obtained solution (from the solve phase) and a local explicit re-evaluation of the PDE at hand using a higher order approximation. Based on the IMEX error indicator, the modified Texas Three Step marking strategy is used to mark the computational nodes for \h-, \p- or \hp-(de-)refinement. Finally, in the refine phase, nodes are repositioned and the order of the method is locally redefined using the variable order of the augmenting monomials according to the instructions from the mark phase.

    The performance of the introduced \hp-adaptive method is first investigated on a two-dimensional Peak problem and further applied to two- and three-dimensional contact problems. We show that the proposed IMEX error indicator adequately captures the global behaviour of the error in all cases considered and that the proposed \hp-adaptive solution procedure significantly outperforms the non-adaptive approach. The proposed \hp-adaptive method stands for another important step towards a fully autonomous numerical method capable of solving complex problems in realistic geometries without the need for user intervention.
}




\maketitle

\section{Introduction}
\label{sec:Introduction}
Many natural and technological phenomena are modelled through Partial Differential Equations (PDEs), which can rarely be solved analytically -- either because of geometric complexity or because of the complexity of the model at hand. Instead, realistic simulations are performed numerically. There are well-developed numerical methods that can be implemented in a more or less effective numerical solution procedure and executed on modern computers to perform virtual experiments or simulate the evolution of various natural or technological phenomena. Nonetheless, despite the immense computing power at our disposal, which allows us to solve ever more complex problems numerically, the development of efficient numerical approaches is still crucial. Relying solely on brute force computing often leads to unnecessarily long computations -- not to mention wasted energy.

Most numerical solutions are obtained using mesh-based methods such as the Finite Volume Method (FVM), the Finite Difference Method (FDM), the Boundary Element Method (BEM) or the Finite Element Method (FEM). Modern numerical analysis is dominated by FEM~\cite{upadhyay2021numerical} as it offers a mature and versatile solution approach that includes all types of adaptive solution procedures~\cite{mitchell2014comparison} and well understood error indicators~\cite{SEGETH20101589}. Despite the widespread acceptance of FEM, the meshing of realistic 3D domains, a crucial part of FEM analysis where nodes are structured into polyhedrons covering the entire domain of interest, is still a problem that often requires user assistance or development of domain-specific algorithms~\cite{liu2005introduction}.

In response to the tedious meshing of realistic 3D domains, required by FEM, and the geometric limitations of FDM and FVM, a new class of mesh-free methods~\cite{liu2002mesh} emerged in the 1970s. Mesh-free methods do not require a topological relationship between computational nodes and can therefore operate on scattered nodes, which greatly simplifies the discretisation of the domain~\cite{zienkiewicz2005finite}, regardless of its dimensionality or shape~\cite{van2021fast, shankar2018robust}. Just recently, they have also been promoted to Computer Aided Design (CAD) geometry aware numerical analysis~\cite{jacquemin_smart_2023}. Moreover, the formulation of mesh-free methods is extremely convenient for implementing \h-refinement~\cite{slak2019adaptive}, considering different approximations of partial differential operators in terms of the shape and size of the \emph{stencil}~\cite{davydov2011adaptive, jacquemin_unified_2021} and the local approximation order~\cite{janvcivc2021monomial}. However, they tend to be more computationally intensive as they require larger stencils for stable computations~\cite{bayona2017role, janvcivc2021monomial} and have limited preprocessing capabilities~\cite{belytschko1996meshless}. This may make them less attractive from a computational point of view, but the ability to work with scattered nodes and easily control the approximation order makes them good candidates for many applications in science and industry~\cite{kosec_simulation_2014, maksic_cooling_2019}.

Adaptive solution procedures are essential in problems where the accuracy of the numerical solution varies spatially and are currently subject of intensive studies. Two conceptually different adaptive approaches have been proposed, namely \p-adaptivity or \h-, \r-adaptivity. In \p-adaptivity, the accuracy of the numerical solution is varied by changing the order of approximation, while in \h- and \r-adaptivity, the resolution of the spatial discretisation is adjusted for the same purpose. In the \h-adaptive approach, nodes are added or removed from the domain as needed, while in the \r-adaptive approach the total number of nodes remains constant -- the nodes are only repositioned with respect to the desired accuracy. Ultimately, \h- and \p-adaptivities can be combined to form the so-called \hp-adaptivity~\cite{gui1985h, gui1986h,devloo12recent}, where the accuracy of the solution is controlled with the order of the method and the resolution of the spatial discretisation.

Since the regions where higher accuracy is required are often not known \emph{a priori}, and to eliminate the need for human intervention in the solution procedure, a measure of the quality of the numerical solution, commonly called a posterior error indicator, is a necessary additional step in an adaptive solution procedure~\cite{liu2005introduction}. The most famous error indicator, commonly referred to as the ZZ-type error indicator, was introduced in 1987 by Zienkiewicz and Zhu~\cite{zienkiewicz1987simple} in the context of FEM and it is still an active research topic~\cite{gonzalez2021error}. The ZZ-type error indicator assumes that the error of the numerical solution is related to the difference between the numerical solution and a locally recovered solution. The ZZ-type error indicator has also been employed in the context of mesh-free solutions of elasticity problems using the mesh-free Finite Volume Method~\cite{Ebrahimnejad} in both weak and strong form using the Finite Point Method~\cite{Angulo}. Furthermore, it also served as an inspiration in the context of Radial Basis Function-Generated Finite Difference (RBF-FD) solution to Laplace equation~\cite{oanh2017adaptive}. Moreover, a residual-based class of error indicators~\cite{Sang} has been demonstrated in the elasticity problems using a Discrete Least Squares mesh-free method~\cite{afshar2011node}. Nevertheless, the most intuitive error indicators are based on the physical interpretation of the solution, usually evaluating the first derivative of the field under consideration~\cite{davydov2011adaptive} or calculating the variance of the field values within the support domain~\cite{slak2019adaptive}.

The advent of \hp-adaptive numerical analysis began with FEM in the 1980s~\cite{guo1986hp}. In \hp-FEM, for example, one has the option of splitting an element into a set of smaller elements or increasing its approximation order. This decision is often considered to be the main difficulty in implementing the \hp-adaptive solution procedure and was already studied by Babuška~\cite{guo1986hp} in 1986. Since then, various decision-making strategies, commonly referred to as marking strategies, have been proposed~\cite{mitchell2014comparison, mitchell2016performance}. The early works use a simple Texas Three Step algorithm, originally proposed in the context of BEM~\cite{tinsley1995three}, where the refinement is based on the maximum value of the error indicator. The first true \hp-strategy was presented by Ainsworth~\cite{ainsworth1997aspects} in 1997, since then many others have been proposed~\cite{mitchell2014comparison,mitchell2016performance}. In general, \p- in FEM is more efficient when the solution is smooth. Based on this observation, most authors nowadays use the local Sobolev regularity estimate to choose between the \h- and the \p-refinement~\cite{houston2005note, houston2005note1, eibner2007adaptive} for a given finite element. Moreover, in~\cite{burg2011convergence} local boundary values are solved, while the authors of~\cite{demkowicz2002fully,rachowicz2006fully} use minimisation of the global interpolation error methods.

For mesh-free methods, \h-adaptivity comes naturally with the ability to work with scattered nodes, and as such has been thoroughly studied in the context of several mesh-free methods~\cite{benito2003h,liu2006stabilized,hu2019spatially}. Only recently, the popular Radial Basis Function-generated Finite Differences (RBF-FD)~\cite{tolstykh2003using} have been used in the \h-adaptive solution of elliptic problems~\cite{oanh2017adaptive,oanh2022approach} and linear elasticity problems~\cite{slak2019adaptive,toth2022h}. Researchers have also reported the combination of \h- and \r-adaptivity, which form a so-called \hr-adaptive solution procedure~\cite{fan2019adaptive}. The \p-adaptive method, on the other hand, is still quite unexplored in the mesh-free community. However, the authors of~\cite{mishra2020adaptive} approach the \p-adaptive RBF-FD method in solving Poisson's equation with the idea of varying the order of the augmenting monomials to maintain the global order of convergence over the domain regardless of the potential variations in the spatial discretisation distance. It should also be noted that some authors reported \p-adaptive methods by locally increasing the number of shape functions, changing the interpolation basis functions, or simply increasing the stencil size~\cite{milewski2021higher, albuquerque2021generalized,liszka1996hp}. These approaches are all to some extent \p-adaptive, but not in their true essence. The authors of~\cite{jancic_p_refined} have introduced a \p-refinement with spatially variable local approximation order and come closest to a true \p-adaptive solution procedure on scattered nodes. However, this work lacks an automated marking and refinement strategy for the local approximation order, e.g.\ based on an error indicator. The automated marking and refinement strategies were used with the weak form \h-\p~adaptive clouds~\cite{duarte1996hp}, where the authors use grid-like \h-enrichment to improve the local field description.

In this paper, we present our attempt to implement the \hp-adaptive strong form mesh-free solution procedure using the mesh-free RBF-FD approximation on scattered nodes. Our solution procedure follows a well-established paradigm based on an iterative loop. To estimate the accuracy of the numerical solution, we employ original IMEX error indicator. The marking strategy used in this work is based on the Texas Three Step algorithm~\cite{eibner2007adaptive}, where the basic idea is to estimate the smoothness or analyticity of the numerical solution. Our refinement strategy is based on the recommendations of~\cite{slak2019adaptive}, where the authors were able to obtain satisfactory results using a purely \h-adaptive solution procedure for elasticity problems. Although the chosen refinement and marking strategies are not optimal~\cite{demkowicz2002fully}, the obtained results clearly outperform the non-adaptive approach.


\section{\hp-adaptive solution procedure}
\label{sec:solution_procedure}
In the present work, we focus on the implementation of mesh-free \hp-adaptive refinement, which combines the advantages of \h- and \p-refinement procedures.
The proposed \hp-adaptive solution procedure follows the well-established paradigm based on an iterative loop, where each iteration step consists of four modules:
\begin{enumerate}
    \item \textbf{Solve} -- A numerical solution $\widehat{u}$ is obtained.
    \item \textbf{Estimate} -- An error indication of the obtained numerical solution.
    \item \textbf{Mark} -- Marking of nodes for refinement/de-refinement.
    \item \textbf{Refine} -- Refinement/de-refinement of the spatial discretisation and local approximation order of the numerical method.
\end{enumerate}

\begin{figure}
    \centering
    \includegraphics[width=0.9\textwidth]{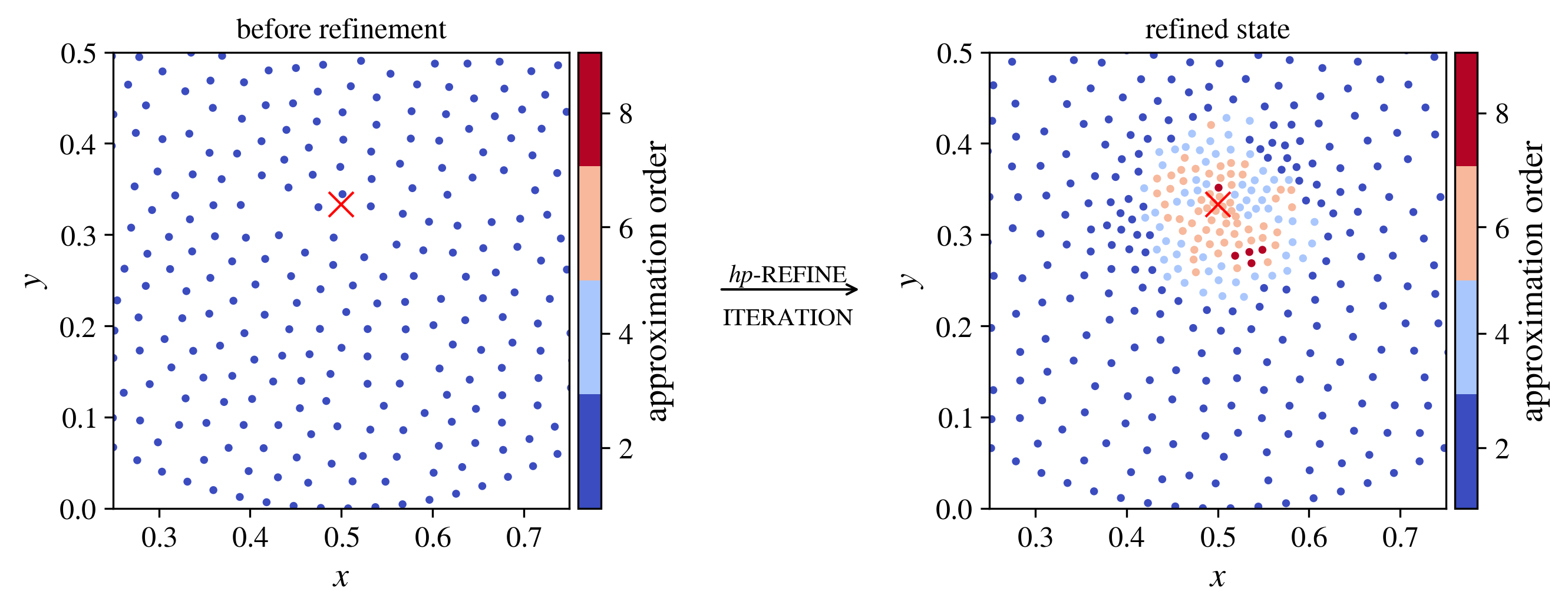}
    \caption{A sketch of a single \hp-refinement iteration for a two-dimensional problem. Note that the exponentially strong source (marked with red cross) is set at $\b p = (\frac{1}{2}, \frac{1}{3})$. The refined state has been obtained by employing \h- and \p-refinement strategies, thus the number of nodes and the local approximation orders in the neighbourhood of the strong source have been modified. Closed form solution has been used to indicate the error in the estimate module.}
    \label{fig:refinement_workflow}
\end{figure}

The workings of each module are further explained in the following subsections, while a full \hp-adaptive solution procedure algorithm is given in Algorithm~\ref{alg:adapt}. For clarity, Figure~\ref{fig:refinement_workflow} also graphically sketches the ultimate goal of a single refinement iteration.
\begin{algorithm}[h]
    \caption{\hp-adaptive solution procedure}
    \label{alg:adapt}
    \vspace{1mm}
    \textbf{Input:} The problem, computational domain $\Omega$, initial nodal density function $h:~\Omega\to \R$, initial approximation order distribution $m:~\Omega\to \N$, the maximal number of iterations $I_{\text{max}}$ and adaptivity parameters $\alpha_{h,p}, \beta_{h,p}, \lambda_{h,p}, \vartheta_{h,p}$. \\
    \textbf{Output:} The \hp-refined numerical solution of the problem.
    \begin{algorithmic}[1]
        \Function{adaptive\_solve}{problem, $\Omega, h, m, I_{\text{max}},\alpha_{h,p}, \beta_{h,p}, \lambda_{h,p}, \vartheta_{h,p}$}
        \For{$i$}{0}{$I_{\text{max}}$} \label{alg:while}
        \State $\Omega^\star \gets \Call{discretise}{\Omega, h}$
        \Comment{Discretises domain using nodal density function $h$.}
        \State $\text{solution} \gets \Call{solve}{\text{problem}, \Omega^\star, m}$
        \Comment{Obtains a numerical solution to the problem.}
        \State $\text{indicator} \gets \Call{imex}{\text{problem}, \text{solution},\Omega^\star, m}$
        \Comment{Error indicator computation.}
        \If{$\Call{stopping\_criteria}{} $}
        \State \Return $\text{solution}$
        \EndIf
        \State $ h, m \gets \Call{adapt}{\text{indicator}, h, m, \Omega^\star, \alpha_{h,p}, \beta_{h,p}, \lambda_{h,p}, \vartheta_{h,p}}$
        \Comment{Refine the nodes and approximation orders.}
        \EndFor
        \State \Return $\text{solution}$
        \EndFunction
    \end{algorithmic}
\end{algorithm}
\subsection{The SOLVE module}
First, a numerical solution $\widehat{u}$ to the governing problem must be obtained. In general, the numerical treatment of a system of PDEs is done in several steps. First, the domain is discretised by positioning the nodes, then the linear differential operators in each computational node are approximated, and finally the system of PDEs is discretised and assembled into a sparse linear system. To obtain a numerical solution $\widehat{u}$, the sparse system is solved.
\paragraph{Domain discretisation}
\begin{wrapfigure}{r}{0.5\textwidth}
    \centering
    \includegraphics[width=0.5\textwidth]{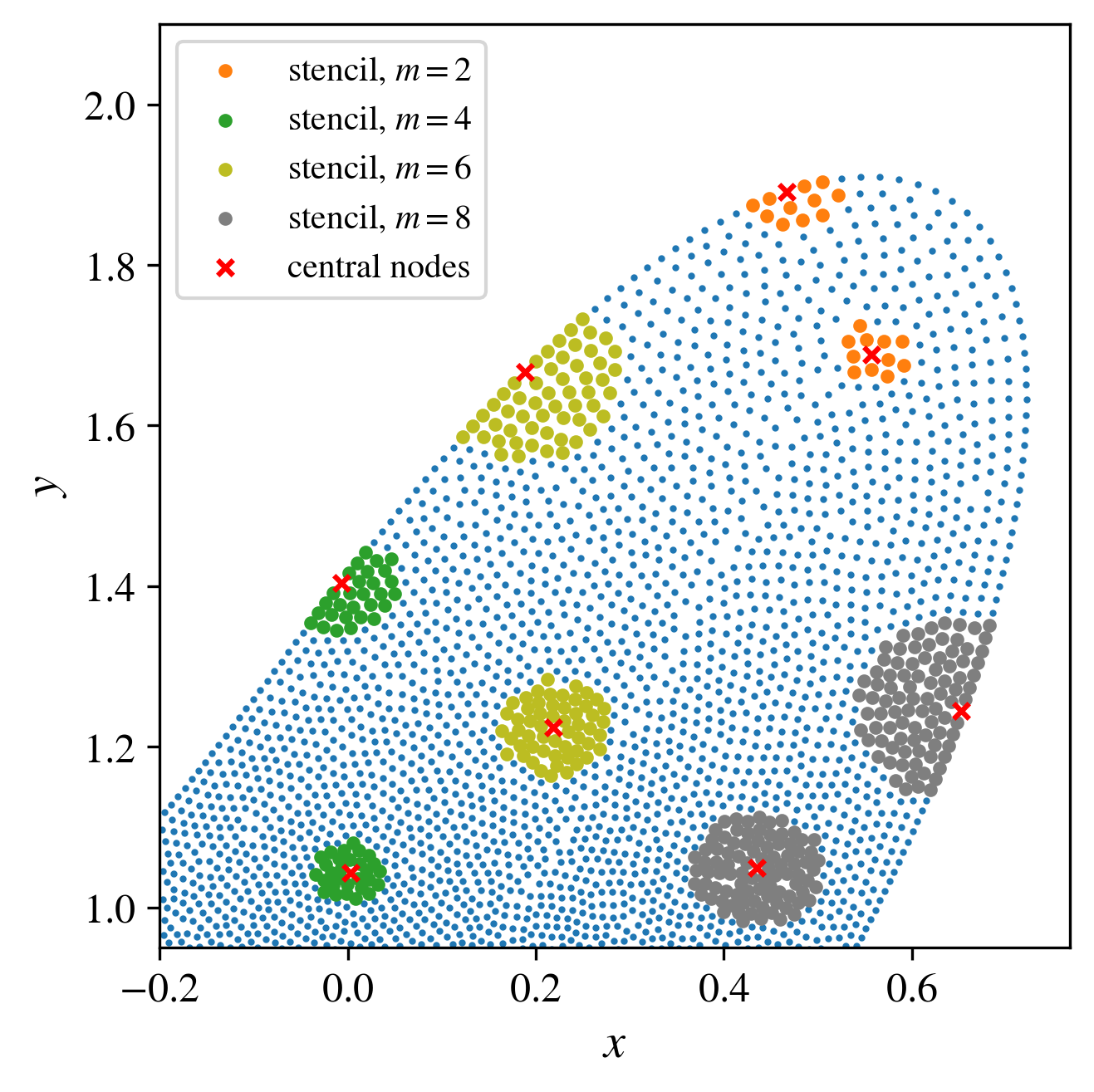}
    \caption{An example of domain discretisation with scattered nodes and variable node density. Example stencils are also shown for different approximation orders $m$ on the domain boundary and its interior.}
    \label{fig:domain}
\end{wrapfigure}
While traditional mesh-based methods discretise the domain by building a mesh, mesh-free methods simplify this step to the positioning of nodes, as no information about internodal connectivity is required. With the mathematical formulation of the mesh-free methods being dimension-independent, we accordingly choose a dimension-independent algorithm for node generation based on Poisson disc sampling~\cite{slak2019generation}. 
Conveniently, the algorithm also supports spatially variable nodal densities required by the \h-adaptive refinement methods. An example of a variable node density discretisation can be found in Figure~\ref{fig:domain}.

Interested readers are further referred to the original paper~\cite{slak2019generation} for more details on the node generation algorithm, its stand-alone C++ implementation in the \emph{Medusa library}~\cite{slak2021medusa}, and follow-up research focusing on its parallel implementation~\cite{depolli_parallel_2022} and parametric surface discretisations~\cite{duh2021fast}.
\paragraph{Approximation of linear differential operators}
Having discretised the domain, we proceed to the approximation of linear differential operators. In this step, a linear differential operator $\L$ is approximated over a set of neighbouring nodes, commonly referred to as \emph{stencil nodes}.

To derive the approximation, we assume a central point $\x _c \in \Omega$ and its stencil nodes $\left \{ \x _i \right\}_{i=1}^n = \mathcal{N}$ for stencil size $n$. A linear differential operator in $\x_c$ is then approximated over its stencil with the following expression
\begin{equation}
    \label{eq:ansatz}
    (\L u)(\x _c) \approx \sum_{i=1}^n w_iu(\x _i),
\end{equation}
for an arbitrary function $u$ and yet to be determined weights $\w$ which are computed by enforcing the equality of approximation~\eqref{eq:ansatz} for a chosen set of basis functions.

In this work, we use Radial Basis Functions (RBFs) augmented with monomials. To eliminate the dependency on a shape parameter, we choose Polyharmonic Splines (PHS) ~\cite{bayona2017role} defined as
\begin{equation}
    f(r) = \begin{cases}r^k,       & k \text{ odd}  \\
             r^k\log r, & k \text{ even}\end{cases},
\end{equation}
for Eucledian distance $r$. The chosen approximation basis effectively results in what is commonly called the RBF-FD approximation method~\cite{tolstykh2003using}.

Furthermore, it is necessary that the stencil nodes form a so-called polynomial unisolvent set~\cite{wendland2004scattered}. In this work, we follow the recommendations of Bayona~\cite{bayona2017role} and define the stencil size as twice the number of augmenting monomials, i.e.
\begin{equation}
    \label{eq:support_size}
    n=2\binom{m+d}{m}
\end{equation}
for monomial order $m$ and domain dimensionality $d$. This, in practice, results in large enough stencil sizes to satisfy the requirement, so that no special treatment was needed to assure unisolvency. While special stencil selection strategies showed promising results~\cite{davydov2011adaptive,DAVYDOV2023115031}, a common choice for selecting a set of stencil nodes $\mathcal{N}$ is to simply select the nearest $n$ nodes. The latter approach was also used in this work. Figure~\ref{fig:domain} shows example stencils for different approximation orders $m$ on domain boundary and its interior.
It is important to note that the augmenting monomials allow us to directly control the order of the local approximation method. The approximation order corresponds to the highest augmenting monomial order $m$ in the approximation basis. However, the greater the approximation order the greater the computational complexity due to larger stencil sizes~\cite{janvcivc2021monomial}. Nevertheless, the ability to control the local order of the approximation method sets the foundation for the \p-adaptive refinement.

To conclude the solve module, the PDEs of the governing problem are discretised and assembled into a global sparse system. The solution of the assembled system stands for the numerical solution $\widehat{u}$.

\subsection{The ESTIMATE module (Implicit-Explicit error indicator)}

In the estimation step, critical areas with high error of the numerical solution are identified. Identifying such areas is not a trivial task. In rare cases where a closed form solution to the governing problem exists, we can directly determine the accuracy of the numerical solution. Therefore, other objective metrics, commonly referred to as \emph{error indicators}, are needed to indicate areas with high error of the numerical solution.
\paragraph{IMplicit-EXplicit (IMEX) error indicator}
In this work we will use an error indicator based on the implicit-explicit~\cite{imex} evaluation of the considered field. IMEX makes use of the implicitly obtained numerical solution and explicit operators (approximated by a higher order basis) to reconstruct the right-hand side of the governing problem.
To explain the basic idea of IMEX, let us define a PDE of type
\begin{equation}
    \label{eq: general PDE}
    \L u = f_{RHS},
\end{equation}
where $\L$ is a differential operator applied to the scalar field $u$ and $f_{RHS}$ is a scalar function. To obtain an error indicator field $\eta$, the problem~\eqref{eq: general PDE} is first solved implicitly by using a lower order approximation $\L^{im}$ of operators $\L$, obtaining the solution $u^{im}$ in the process. The explicit high order operators $\L^{ex}$ are then used over the implicitly computed field $u^{im}$ to reconstruct the right-hand side of the problem~\eqref{eq: general PDE} obtaining $f_{RHS}^{ex}$ in the process.
The error indication is then calculated as $\eta = \lvert f_{RHS} - f_{RHS}^{ex} \rvert$.
The calculation steps of the IMEX error indicator are also shown in Algorithm~\ref{alg:imex}.
\begin{algorithm}[h]
    \caption{IMEX error indicator}
    \label{alg:imex}
    \vspace{1mm}
    \textbf{Input:} The problem, domain $\Omega$, differential operators $\L$, low-order approximation basis $\xi$, high order approximation basis $\zeta$.\\
    \textbf{Output:} Error indicator field $\eta$.
    \begin{algorithmic}[1]
        \Function{indicate\_error}{problem, $\Omega$, $\L$, $\xi$, $\zeta$}
        \State $\L^{im} \gets \Call{approximate}{\Omega, \xi}$
        \Comment{Obtain low-order approximation of differential operators $\L$.}
        \State $u^{im}, f_{RHS} \gets \Call{solve}{\text{problem}, \Omega, \L^{im}}$
        \Comment{Obtain a numerical solution to the problem.}
        \State $\L^{ex} \gets \Call{approximate}{\Omega, \zeta}$
        \Comment{Obtain high order approximation of differential operators $\L$.}
        \State $f_{RHS}^{ex} \gets \Call{evaluate}{\text{problem}, \Omega, \L^{ex}, u^{im}}$
        \Comment{Explicit re-evaluation.}
        \State $ \eta \gets \Call{compute}{f_{RHS}, f_{RHS}^{ex}}$
        \Comment{Obtain error indicator field.}
        \State \Return $\eta$
        \EndFunction
    \end{algorithmic}
\end{algorithm}

The assumption that the deviation of the explicit high order evaluation $\L^{ex} u^{im}$ from the exact $f_{RHS}$ corresponds to the error of the solution $u^{im}$ is similar to the reasoning behind the ZZ-type indicators, where the deviation of the recovered high order solution from the computed solution characterises the error. As long as the error in $u^{im}$ is high, the explicit re-evaluation will not correctly solve the Equation~\eqref{eq: general PDE}. However, as the error in $u^{im}$ decreases, the difference between $f_{RHS}$ and $f_{RHS}^{ex}$ will also decrease, assuming that the error is dominated by the inaccuracy of $u_{im}$ and not by the differential operator approximation.

It is worth noting that the definition of IMEX is general in the sense that computing the error indication $\eta$ does not distinguish between the interior and boundary nodes. In the boundary nodes, the error indicator $\eta$ is calculated in the same way as in the interior nodes. In the case of Dirichlet boundary conditions, the error indicator is trivial because the solution fields are exactly imposed, i.e.\ the error indicator results in $\eta = 0$. However, in case of boundary conditions involving the evaluation of derivatives (Robin and Neumann), $\eta \neq 0$.


\subsection{The MARK module}
\label{sec:mark}
After the error indicator $\eta$ has been obtained for each computational point in domain $\Omega$, a marking strategy is applied. The main goal of this module is to mark the nodes with too high or too low values of the error indicator in order to achieve a uniformly distributed accuracy of the numerical solution and to reduce the computational cost of the solution procedure -- by avoiding fine local field descriptions and high order approximations where this is not required. Moreover, the marking strategy not only decides whether or not (de-)refinement should take place at a particular computational node, but also defines the type of refinement procedure if there are several to choose from. In this work, we use a modified Texas Three Step marking strategy~\cite{tinsley1995three,heuer2001hp}, originally restricted to refinement (no de-refinement) with the \h- and \p-refinement types. This chosen strategy was also considered in one of the recent papers by Eibner~\cite{eibner2007adaptive}, who showed that, although extremely simple to understand and implement, it can provide results good enough to demonstrate the advantages of mesh-based \hp-adaptive solution procedures.

In each iteration of the adaptive procedure, the marking strategy starts by checking the error indicator values $\eta _i$ for all computational nodes in the domain. Unlike the originally proposed marking strategy~\cite{eibner2007adaptive} that used only refinement, we additionally introduce de-refinement. Therefore, if $\eta _i$ is greater than $\alpha \eta_{max}$ for the maximum indicator value $\eta_{max}$ and a free model parameter $\alpha \in (0, 1)$, the node is marked for refinement. If $\eta _i$ is less than $\beta \eta_{max}$ for a free model parameter $\beta \in (0,1) \land \beta \leq \alpha$, the node is marked for de-refinement. Otherwise, the node remains unmarked, which means that no (de-)refinement should take place. The marking strategy can be summarised with a single equation
\begin{equation}
    \label{eq:marking_strategy}
    \begin{cases}
        \eta _i > \alpha \eta_{max},                          & \text{ refine}     \\
        \beta \eta_{max} \leq \eta _i \leq \alpha \eta_{max}, & \text{ do nothing} \\
        \eta _i < \beta \eta_{max},                           & \text{ de-refine}
    \end{cases}.
\end{equation}

In the context of mesh-based methods, it has already been observed, that such marking strategy, although easy to implement, is far from optimal~\cite{eibner2007adaptive,mitchell2014comparison}. Additionally, it has also been demonstrated that in case of smooth solutions \p-refinement is preferred while  \h-refinement is preferred in volatile fields, e.g.\ in vicinity of a singularity in the solution~\cite{mitchell2014comparison,demkowicz2002fully}, which cannot be achieved with the chosen marking strategy. 
Additional discussion on this issue can be found in Section~\ref{sec:elasticity}, where problems with singularity in the solution are discussed, and in Section~\ref{sec:step_beyond} where we discuss some guidelines for possible work on improved marking strategies.

\begin{figure}
    \centering
    \includegraphics[width=0.9\textwidth]{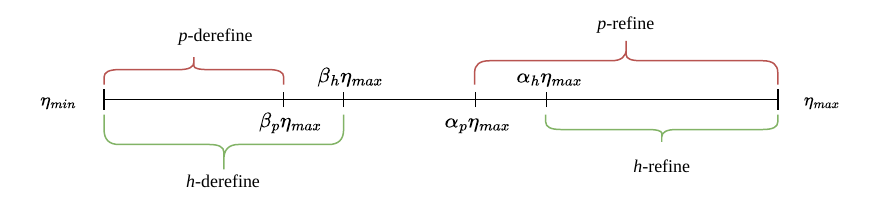}
    \caption{A visual representation of \h- and \p-(de)refinement marking strategy.}
    \label{fig:indicator_complex}
\end{figure}

Since our work is focused on the implementation of \hp-adaptive solution procedure rather than discussing the optimal marking strategy, we decided to secure full control over the marking strategy by treating \h- and \p-methods separately -- but at the cost of higher number of free parameters. Therefore, the marking strategy is modified by introducing parameters $\left \{ \alpha_h, \beta_h \right \}$ and $\left \{ \alpha_p, \beta_p \right \}$ for separate treatment of \h- and \p-refinements, respectively (see Figure~\ref{fig:indicator_complex} for clarification). Note that the proposed modified marking strategy can mark a particular node for \h-, \p- or \hp-(de-)refinement if required, otherwise the computational node is left unchanged.


\subsection{The REFINE module}
\begin{figure}
    \centering
    \includegraphics[width=0.9\textwidth]{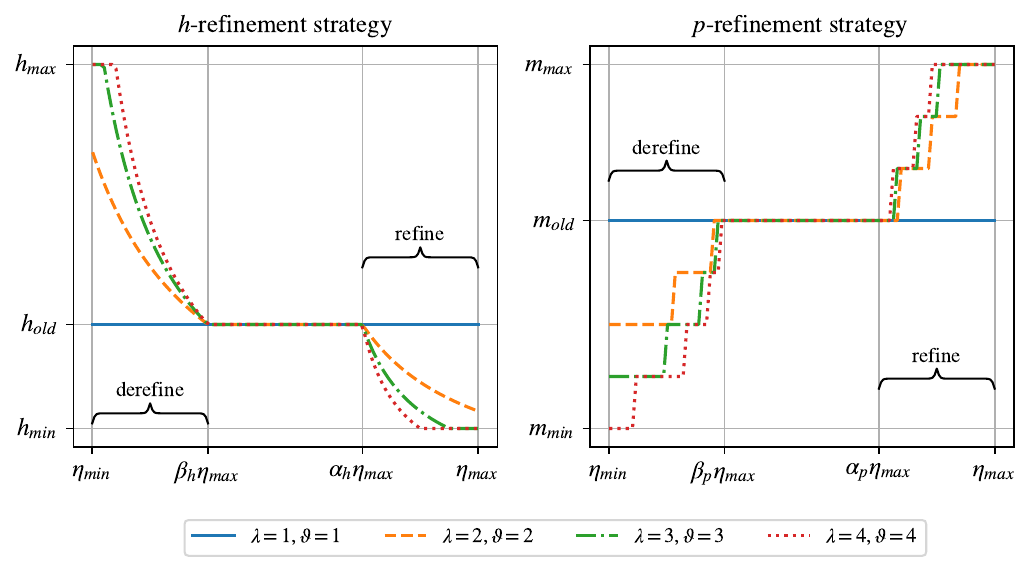}
    \caption{A visual representation of the (de-)refinement strategies for different values of refinement aggressiveness $\lambda$ and de-refinement aggressiveness $\vartheta$. Notice that both refinement types also have lower and upper limits.}
    \label{fig:refinement_rules}
\end{figure}

After obtaining the list of nodes marked for modification, the refinement module is initialised. In this module, the local field description and local approximation order are left unchanged for the unmarked nodes, while the remaining nodes are further processed to determine other refinement-type-specific details -- such as the amount of the (de-)refinement. Our $h$-refinement strategy is inspired by the recent $h$-adaptive mesh-free solution of elasticity problem~\cite{slak2019adaptive}, where the following \h-refinement rule was introduced
\begin{equation}
    \label{eq:refinement}
    h_i^{new}(\b p) = \frac{h_i^{old}}{\frac{\eta_i - \alpha \eta _{max}}{\eta_{max} - \alpha \eta_{max}}\Big(\lambda - 1\Big) + 1}
\end{equation}
for the dimensionless parameter $\lambda \in [1, \infty)$ allowing us to control the aggressiveness of the refinement -- the larger the value, the greater the change in nodal density, as shown in Figure~\ref{fig:refinement_rules} on the left. This refinement rule also conveniently refines the areas with higher error indicator values more than those closer to the upper refinement threshold $\alpha_h \eta_{max}$. Similarly, a de-refinement rule is proposed
\begin{equation}
    \label{eq:derefinement}
    h_i^{new}(\b p) = \frac{h_i^{old}}{\frac{\beta \eta _{max} - \eta_i}{\beta\eta_{max} - \eta_{min}}\Big(\frac{1}{\vartheta} - 1\Big) + 1},
\end{equation}
where parameter $\vartheta \in [1, \infty)$ allows us to control the aggressiveness of de-refinement.

The same refinement~\eqref{eq:refinement} and de-refinement~\eqref{eq:derefinement} rules are applied to control the order of local approximation ($p$-refinement), except that this time the value is rounded to the nearest integer, as shown in Figure~\ref{fig:refinement_rules} on the right. Similarly, and for the same reasons as for the marking strategy (see Section~\ref{sec:mark}), we consider a separate treatment of \h- and \p-adaptive procedures by introducing (de-)refinement aggressiveness parameters $\left \{\lambda_h, \vartheta_h \right \}$ and $\left \{\lambda_p, \vartheta_p \right \}$ for \h- and \p-refinement types respectively.

\subsection{Finalization step}
Before the 4 modules can be iteratively repeated, the domain is re-discretised taking into account the newly computed local internodal distances $h_i^{new}(\b p)$ and the local approximation orders $m_i^{new}(\b p)$. However, both are only known in the computational nodes, while global functions $\widehat{h}^{new}(\b p)$ and $\widehat{m}^{new}(\b p)$ over our entire domain space $\Omega$ are required.

We use Sheppard's inverse distance weighting interpolation using the closest $n^h_s$ neighbours to construct $\widehat{h}^{new}(\b p)$ and the closest $n^m_s$ neighbours to construct $\widehat{m}^{new}(\b p)$. In general, the proposed refinement strategy can introduce aggressive and undesirable local jumps in node density, which ultimately leads to a potential violation of the quasi-uniform internodal spacing requirement within the stencil. To mitigate this effect, we use relatively large $n^h_s = 30$ to smoothen such potential local jumps. The $\widehat{m}^{new}(\b p)$ is much less sensitive in this respect and therefore a minimum $n^m_s = 3$ is used.

\begin{figure}
    \centering
    \includegraphics[width=0.95\textwidth]{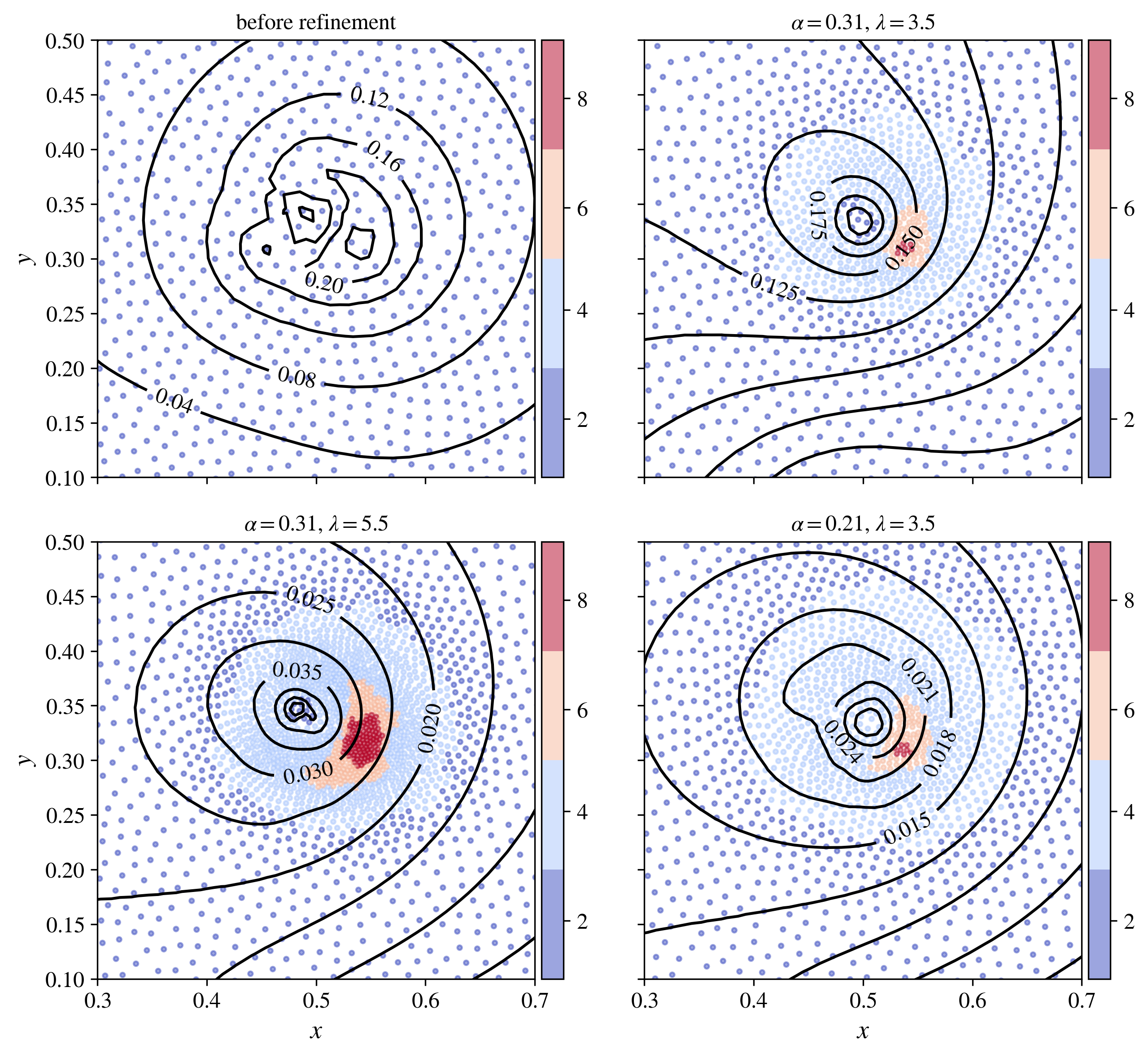}
    \caption{Demonstration of \hp-refinement for selected values of refinement parameters. The top left figure shows the numerical solution before its refinement, while the rest show its refined state for different values of refinement parameters. Contour lines are used to show the absolute error of the numerical solution. To denote the \p-refinement, the nodes are coloured according to the local approximation order. For clarity, all figures are zoomed to show only the neighbourhood of an exponentially strong source $e^{-a \left \| \x - \x_s \right \|^2}$ positioned at $\x_s = \Big(\frac{1}{2}, \frac{1}{3}\Big)$.}
    \label{fig:refinement_demonstration}
\end{figure}

Figure~\ref{fig:refinement_demonstration} schematically demonstrates 3 examples of \hp-refinements. For demonstration purposes, the refinement parameters for \h- and \p-adaptivity are the same, i.e.\ $\left \{ \alpha, \beta, \lambda, \vartheta \right \} = \left \{ \alpha_h, \beta_h, \lambda_h, \vartheta_h \right \} = \left \{ \alpha_p, \beta_p, \lambda_p, \vartheta_p \right \}$. Additionally, the de-refinement aggressiveness $\vartheta$ and the lower threshold $\beta$ are kept constant, so that effectively only the upper limit of refinement $\alpha$ and the refinement aggressiveness $\lambda$ are altered. We observe that the effect of the refinement parameters is somewhat intuitive. The greater the aggressiveness $\lambda$, the better the local field description and the greater the number of nodes with high approximation order. A similar effect is observed when manipulating the upper refinement threshold $\alpha$, except that the effect comes at a smoother manner. Note also that all refined states were able to increase the accuracy of the numerical solution from the initial state.

\subsection{Note on marking and refinement strategies}
\label{sec:note_mark_refine}
With the chosen marking and refinement strategies, a separate treatment of \h- and \p-refinement types turned out to be a necessary complication for a better overall performance of the solution procedure. Nevertheless, we have tried to simplify the solution procedure as much as possible. In the process, important observations have been made -- some of which we believe should be highlighted. This section therefore opens a discussion on important remarks related to the proposed marking and refinement modules.

\subsubsection{The error indicators}
Since the \h- and \p-refinements are conceptually different, our first attempt was to employ two different error indicators -- one for each type of refinement. We employed the previously proposed variance of field values~\cite{slak2019adaptive} for marking the \h-refinement and the approximation order based IMEX for the \p-refinement. Unfortunately, no notable advantages of such solution procedure has been observed and was therefore discarded due to the increased implementation complexity. However, other combinations that might show more promising results should be considered in future work.

\subsubsection{Free parameters}
In the proposed solution procedure, each adaptivity type comes with $4$ free parameters that need to be defined, i.e.\ $\left \{ \alpha_{h,p}, \beta_{h,p}, \lambda_{h,p}, \vartheta_{h,p} \right \}$. This gives a total of $8$ free parameters that can be fine-tuned to a particular problem. While we have tried to avoid any kind of fine-tuning, we have nevertheless observed that these parameters can have a crucial impact on the overall performance of the \hp-adaptive solution procedure in terms of (i) the achieved accuracy of the numerical solution, (ii) the spatial variability of the error of the numerical solution, (iii) the computational complexity, and (iv) the stability of the solution procedure.

We observed that if the refinement aggressiveness $\lambda_h$ is too high, the number of nodes can either diverge into unreasonably large domain discretisations or ultimately violate the quasi-uniform internodal spacing requirement, making the solution procedure unstable. Note that here we refer to the stability of the solution of the discretized PDEs, which ultimately governs the stability of the whole solution procedure.
Furthermore, a large number of nodes combined with high approximation orders can lead to unreasonably high computational complexity in a matter of few iterations. However, when refinement aggressiveness $\lambda_h$ and $\lambda_p$ is set too low, the number of required iterations can increase to such an extent that the entire solution procedure becomes inefficient. On top of that, the lower and upper threshold multipliers $\alpha$ and $\beta$ also play a crucial role. If $\alpha$ is too low, almost the entire domain is refined.
Moreover, if $\alpha$ is too large, almost no refinement takes place and if it does, it is extremely local, which again has no beneficial consequences as it often leads to a violation of the quasi-uniform nodal distribution requirement.

In our tests, based on extensive experimental parameter testing, we have selected a reasonable combination of all $8$ parameters that lead to a stable solution procedure while demonstrating the advantages of the proposed \hp-adaptive approach.
A thorough analysis of these parameters and their correlation would most likely lead to better results, as there is no guarantee that the selected parameters are optimal. However, such an analysis is beyond the scope of this paper, whose aim is to present an \hp-adaptive solution procedure in the context of mesh-free methods and not to discuss the optimal marking and refinement strategies.
Nevertheless, we have tried to reduce the number of free parameters by using the same values for \h- and \p-adaptivity (see Figure~\ref{fig:refinement_demonstration}). While this approach also yielded satisfactory results that outperformed the numerical solutions obtained with uniform nodal and approximation order distributions in terms of accuracy, the full 8-parameter formulation easily yielded significantly better results.
%
\subsubsection{A step beyond the artificial refinement strategies}
\label{sec:step_beyond}
As discussed in Section~\ref{sec:mark} and later in Section~\ref{sec:elasticity}, the Texas Three Step based marking strategy cannot assure the optimal balance of \h- and \p-refinements due to missing local data regularity estimation~\cite{mitchell2014comparison}. In FEM, local Sobolev regularity estimate is commonly used to choose between the \h- and the \p-refinement~\cite{houston2005note, houston2005note1, eibner2007adaptive}. Using an estimate for upper error bound~\cite{bayona2019insight,tominec} one could generalise this approach to meshless methods, essentially upgrading the strategy with an information on the minimal internodal spacing required for local approximation of the partial differential operator of a certain order.

The refinement strategy could also be based on a specific knowledge about convergence rates and computational complexity in terms of internodal distance $h(\b p)$ and local approximation orders $m(\b p)$.

It has already been shown by Bayona~\cite{bayona2019comparison} that the approximation error of mesh-free interpolant $F$ is bounded by
\begin{equation}
    \label{eq:errr}
    \left \| F(\b p) - u(\b p) \right \| _\infty \leq C h^{m +1}\max _{\b p \in \Omega}\lvert \L^{(m  +1 )}(u(\b p)) \rvert .
\end{equation}
Note that the constant $C$ present in Equation~\eqref{eq:errr} depends on the stencil and on the approximation order, both of which are modified by the \hp-adaptive solution procedure. Nevertheless, for the purpose of illustrating how a better marking strategy could be constructed, we decide to simplify the Equation~\eqref{eq:errr} to saying that the error $e$ is proportional to $h(\b p)^{m(\b p)}$. Knowing the target error $e_t$, we write the ratio of $e_t/e_0$ as
\begin{equation}
    \label{eq:p_guess}
    \frac{e_t}{e_0}\propto \frac{h^{m_t}}{h^{m_0}} = h^{m_t -m_0},
\end{equation}
where $m_t$ is used to denote the target approximation order and $m_0$ is the current order of the approximation used to compute current error $e_0$.

From Equation~\eqref{eq:p_guess} a smarter guess for target local approximation order can be obtained
\begin{equation}
    \label{eq:p_guess_log}
    m_t = m_0 + \ln{\frac{e_t}{e_0}}.
\end{equation}
Such strategy would conveniently leave the approximation order unchanged when $e_t = e_0$, increase it when $e_t < e_0$ and decrease it when $e_t > e_0$.

A step even further could be to additionally consider the change in computational complexity, similar to what the authors of~\cite{burg2011convergence} and~\cite{mishra2020adaptive} have already shown. Therefore, we believe that future work should consider the minimum local computational complexity criteria. A rough computational complexity can be obtained with the help of
\begin{equation}
    \label{eq:complexity}
    \chi \propto \frac{\binom{m_t + d}{d}^3\Big(\frac{1}{h_t}\Big)^d}{\binom{m_0 + d}{d}^3\Big(\frac{1}{h_0}\Big)^d},
\end{equation}
for domain dimensionality $d$ and target and current internodal distances $h_t$ and $h_0$ respectively.

\subsection{Implementation note}
\label{sec:implementation}

The entire \hp-adaptive solution procedure from Algorithm~\ref{alg:adapt} is implemented in C++. All meshless methods and approaches used in this work are included in our in-house developed \emph{Medusa library}~\cite{slak2021medusa}. 
The code\footnote{The source code is available at: \url{https://gitlab.com/e62Lab/public/2022_p_hp-adaptivity} under tag \emph{v1.2}.} was compiled using \texttt{g++ (GCC) 9.3.0 for Linux} with \texttt{-O3 -DNDEBUG -fopenmp} flags. Post-processing was done using Python 3.10 and Jupyter notebooks, also available in the provided git repository.

\section{Demonstration on exponential peak problem}
\label{sec:demonstration}
The proposed \hp-adaptive solution procedure is first demonstrated on a synthetic example. We chose a 2-dimensional Poisson problem with exponentially strong source positioned at $\x_s = \Big (\frac{1}{2}, \frac{1}{3}\Big )$. This example is categorized as a difficult problem and is commonly used to test the performance of adaptive solution procedures~\cite{mitchell2014comparison,oanh2022approach,daniel2018adaptive,mitchell2016performance}. The problem has a tractable solution $u(\x)=e^{-a \left \| \x - \x_s \right \|^2}$, which allows us to evaluate the precision of the numerical solution $\widehat{u}$, e.g.\ in terms of the infinity norm
\begin{equation}
    e_{\infty} = \frac{\|\widehat u - u\|_\infty}{\|u\|_\infty}, \quad \|u\|_\infty = \max_{i=1, \ldots, N} \lvert u_i \rvert.
\end{equation}

\begin{wrapfigure}{r}{0.5\textwidth}
    \centering
    \includegraphics[width=0.5\textwidth]{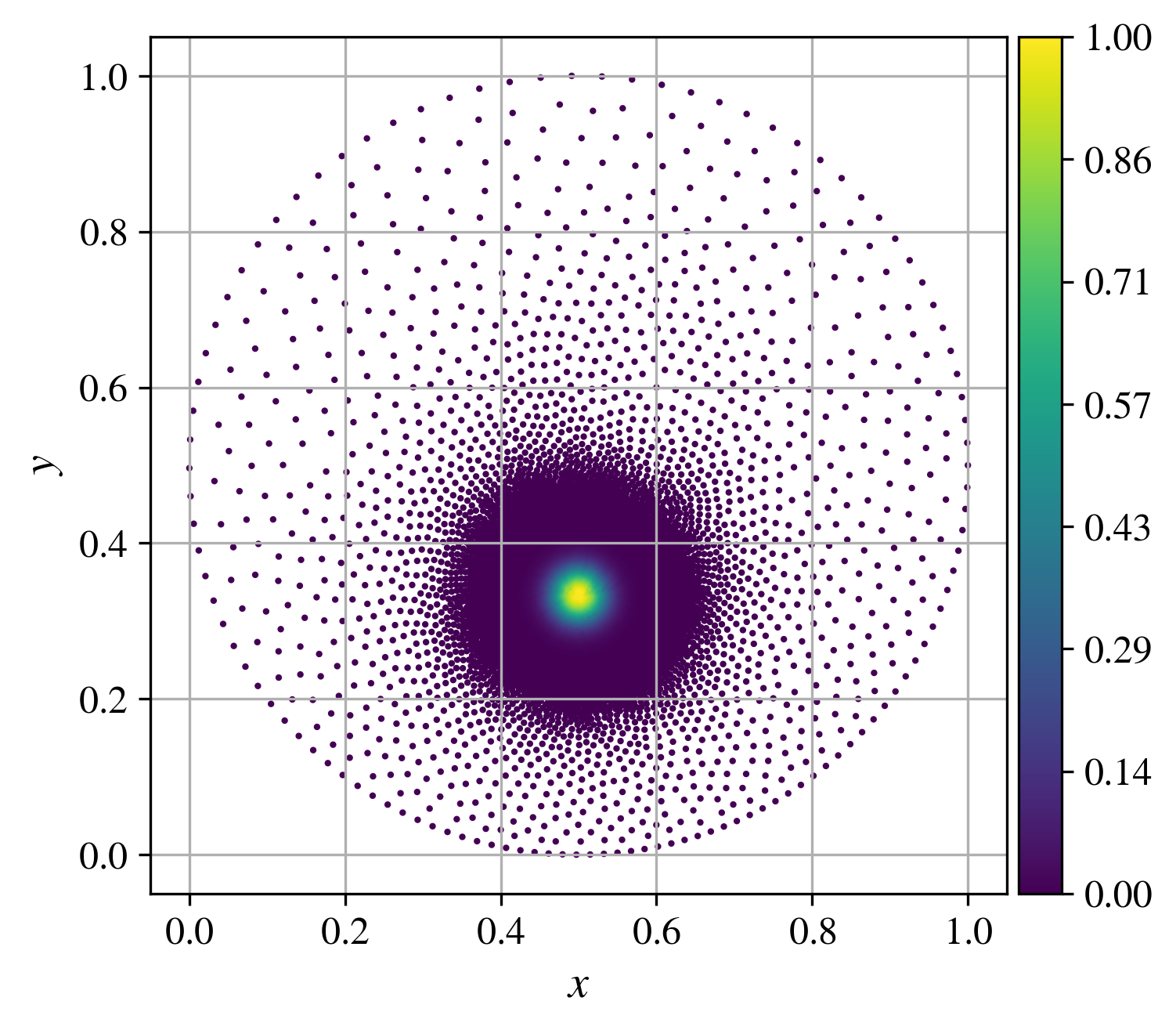}
    \caption{Example \hp-refined solution to exponential peak problem.}
    \label{fig:solution_example}
\end{wrapfigure}

Governing equations are
\begin{align}
    \lap u (\x)   & = 2a e^{-a \left \| \x - \x_s \right \|^2}(2a\left \| \x - \x_s \right \| - d) & \text{in } \Omega,   \\
    u (\x)        & = e^{-a \left \| \x - \x_s \right \|^2}                                        & \text{on } \Gamma_d, \\
    \nabla u (\x) & = -2a(\x - \x_s)e^{-a \left \| \x - \x_s \right \|^2}                          & \text{on } \Gamma_n,
\end{align}
for a $d$-dimensional domain $\Omega$ and strength $a=10^3$ of the exponential source. The domain boundary is split into two sets: Neumann $\Gamma_n =\left \{ \x, x \leq \frac{1}{2} \right \}$ and Dirichlet $\Gamma_d =\left \{ \x, x > \frac{1}{2}\right \}$ boundaries. An example \hp-refined numerical solution is shown in Figure~\ref{fig:solution_example}.

In the continuation of this paper, the numerical solution of the final linear system is obtained by employing BiCGSTAB solver with a ILUT preconditioner from the Eigen C++ library~\cite{eigenweb}. Global tolerance was set to $10^{-15}$ with a maximum number of $800$ iterations and drop-tolerance and fill-factor set to $10^{-5}$ and $50$ respectively. While the initial adaptivity solution was obtained without the guess, all other iterations used the previous numerical solution $\widehat{u_{i-1}}$ as the guess for new numerical solution $\widehat{u_{i}}$, effectively reducing the number of iterations required by the BiCGSTAB solver.

\subsection{Convergence analysis of unrefined solution}
\label{sec:unrefined}

The problem is first solved on a two-dimensional unit disc without employing any refinement procedures, i.e.\ with uniform nodal and approximation order distributions. The shapes approximating the linear differential operators are computed using the RBF-FD with PHS order $k=3$ and monomial augmentation $m\in \left \{2, 4, 6, 8 \right \}$.

\begin{figure}
    \centering
    \includegraphics[width=0.9\textwidth]{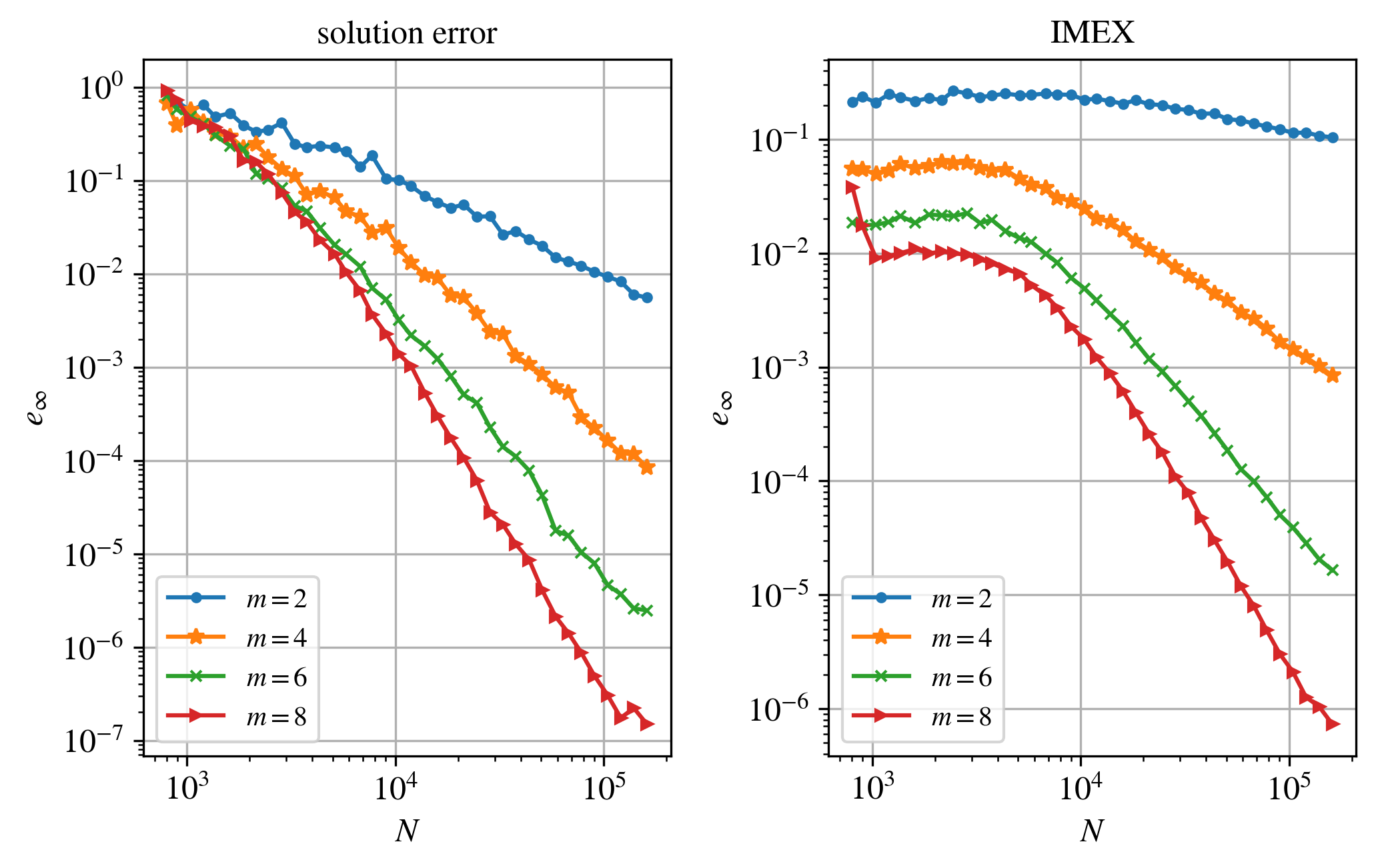}
    \caption{Convergence of unrefined numerical solution (left) and IMEX error indicator (right). Figure only shows a median value after 50 runs with slightly different domain discretisations. Note that, the approximation order $m$ in the right figure denotes the approximation order used to obtain the numerical solution, while the explicit operators employed by the IMEX error indicator are approximated with orders $m+2$.}
    \label{fig:convergence_uniform}
\end{figure}

Figure~\ref{fig:convergence_uniform} shows the results. Each plotted point is an average obtained after 50 consecutive runs with slightly different domain discretisations (a random seed for generating expansion candidates was changed, see~\cite{slak2019generation} for more details). In this way, we can not only study the convergence behaviour, but also evaluate how prone the numerical method is to non-optimal domain discretisations. The convergence of the numerical solution for selected monomial augmentations is shown on the left. We observe that due to the strong source, the convergence rates no longer follow the theoretical prediction of being proportional to $h^m$. Instead, the convergence rates for a small number of computational nodes ($N\lessapprox 2000$) are significantly lower than that obtained for larger domain discretisations ($N\gtrapprox 3000$) for all approximation orders $m > 2$. Furthermore, the accuracy gain by using higher order approximations with small domain discretisations is practically negligible.
However, when the local field description is sufficient, both the numerical solution and the IMEX error indicator (Figure~\ref{fig:convergence_uniform} on the right) give reliable results. While we could have forced at least one node in the neighbourhood of the source, we do not use any special techniques in this work. Instead, further research is simply limited to sufficiently large domains so that this observation does not represent an issue.

Moreover, the behaviour of the IMEX error indicator is studied on the right side of Figure~\ref{fig:convergence_uniform}. Here, the approximation order $m$ means that the implicit numerical solution $u^{im}$ was obtained with approximation order $m$, while the explicit operators $\L^{ex}$ from IMEX were approximated using monomials up to and including order $m+2$. The observations show that the maximum value of the error indicator also converges with the number of computational nodes. Moreover, we can also observe the aforementioned change in the convergence rate of the numerical solution, since the maximum value of the error indicator for domain sizes $N\lessapprox 3000$ is approximately constant.


%
%
\subsection{Analysis of \hp-refined solution}
The same problem is now solved by employing the \hp-adaptivity. Free parameters are adjusted to each refinement type, as can be seen in Table~\ref{tab:params}. Adaptivity iteration loop is stopped after a maximum of $N_{iter}$ iterations. For practical use, other stopping criteria could also be used, e.g.\ based on the maximum error indicator reduction
\begin{equation} \frac{\eta_{max}^j}{\eta_{max}^0} \leq \gamma,
\end{equation}
for the iteration index $j$. The shapes are computed with RBF-FD using the PHS with order $k=3$ and local monomial augmentation restricted to choose between approximation orders $m\in \left \{2, 4, 6, 8 \right \}$. Note that the IMEX error indicator increases the local approximation order by $2$, effectively using monomial orders $m_{IMEX}\in \left \{4, 6, 8, 10 \right \}$. Furthermore, to avoid unreasonably large number of computational nodes, the maximum number of allowed nodes $N_{max}$ is defined. Once this number is reached, further \h-refinement is prevented and only de-refinement is allowed, while the \p-adaptive method retains its full functionality. To avoid insufficient local field description, the local nodal density is limited by an upper bound, i.e.\ $h(\b p ) \leq h_{max}$. The order of the PHS is left constant.

\begin{table}
    \centering
    \renewcommand{\arraystretch}{1.2}
    \begin{tabular}{cccc|cccc||ccc} \hline
        \multicolumn{1}{c}{$\beta_h$} & \multicolumn{1}{c}{$\alpha_h$} & \multicolumn{1}{c}{$\lambda_h$} & \multicolumn{1}{c|}{$\vartheta_h$} & \multicolumn{1}{c}{$\beta_p$} & \multicolumn{1}{c}{$\alpha_p$} & \multicolumn{1}{c}{$\lambda_p$} & \multicolumn{1}{c||}{$\vartheta_p$} & \multicolumn{1}{c}{$h_{max}$} & \multicolumn{1}{c}{$N_{max}$} & \multicolumn{1}{c}{$N_{iter}$} \\ \hline
        0.175                         & 0.225                          & 2.625                           & 1.01                               & $10^{-4}$                     & 0.05                           & 5                               & 1.258                               & 0.1                           & $2.5\cdot 10^{5}$             & 70                             \\ \hline
    \end{tabular}
    \caption{Adaptivity parameters used to obtain solution to the peak problem.}
    \label{tab:params}
\end{table}

\subsubsection{A brief analysis of IMEX error indicator}
\label{sec:imex_analysis}
\begin{figure}
    \centering
    \includegraphics[width=\textwidth]{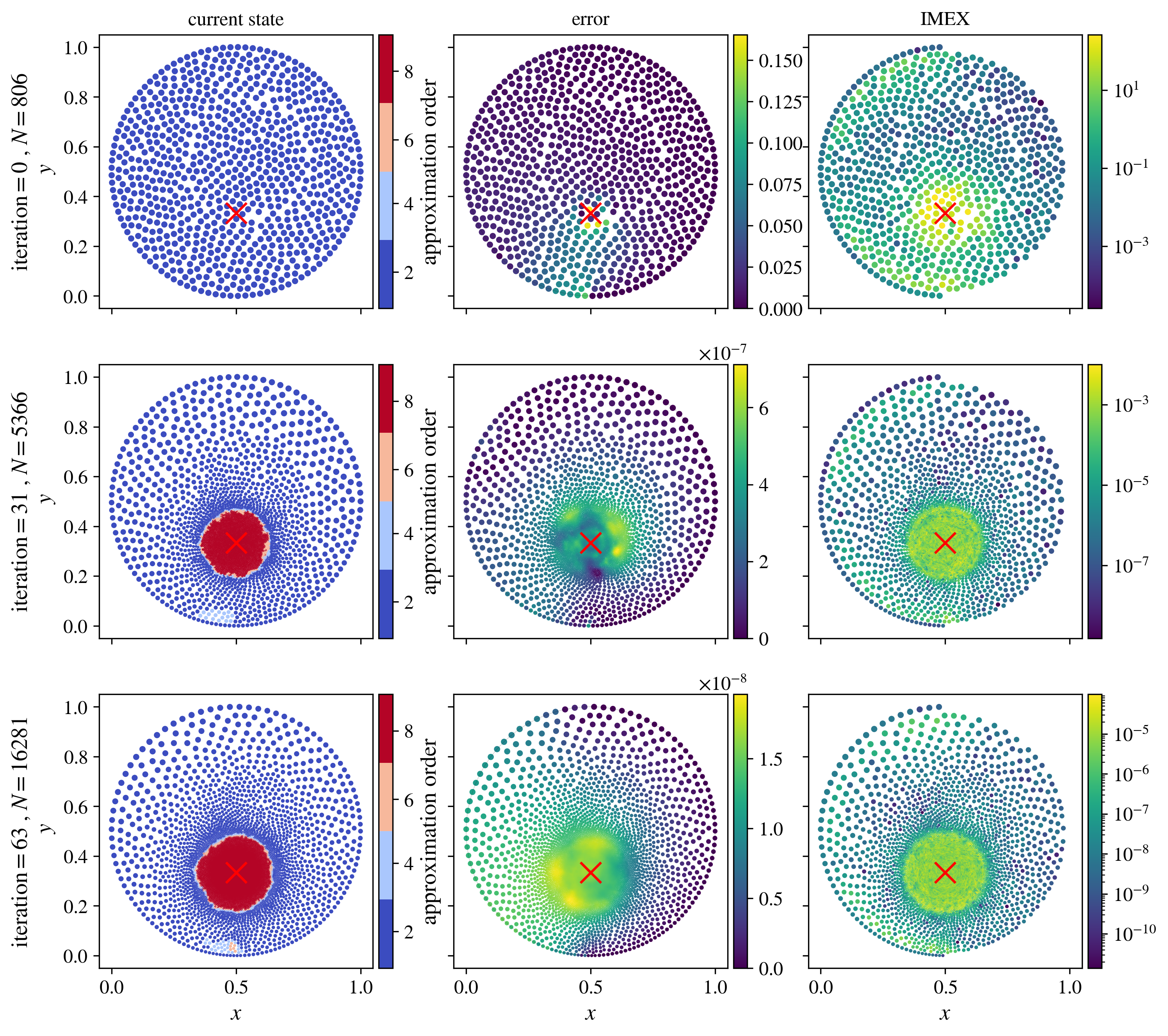}
    \caption{Refinement demonstration. Initial iteration (top row), intermediate iteration (middle row) and best-performing iteration (bottom row) accompanied with solution error (middle column) and IMEX error indicator values (right column). The IMEX values for Dirichlet boundary nodes are not shown. A red cross is used to mark the location of the strong peak.}
    \label{fig:ref_imex_demonstration}
\end{figure}

Figure~\ref{fig:ref_imex_demonstration} shows example indicator fields for the initial iteration, the intermediate iteration, and the iteration that achieved the best numerical solution accuracy -- hereafter also referred to as the best-performing iteration or simply the best iteration. The third column shows the IMEX error indicator. We can see that the IMEX has successfully located the position of the strong source at $\x_s = \Big (\frac{1}{2}, \frac{1}{3}\Big )$ as the highest indicator values are seen in its vicinity. Furthermore, the second column shows that both the accuracy of the numerical solution and the uniformity of the error distribution were significantly improved by the \hp-adaptive solution procedure, further proving that IMEX can be successfully used as a reliable error indicator.

\begin{figure}
    \centering
    \includegraphics[width=0.9\textwidth]{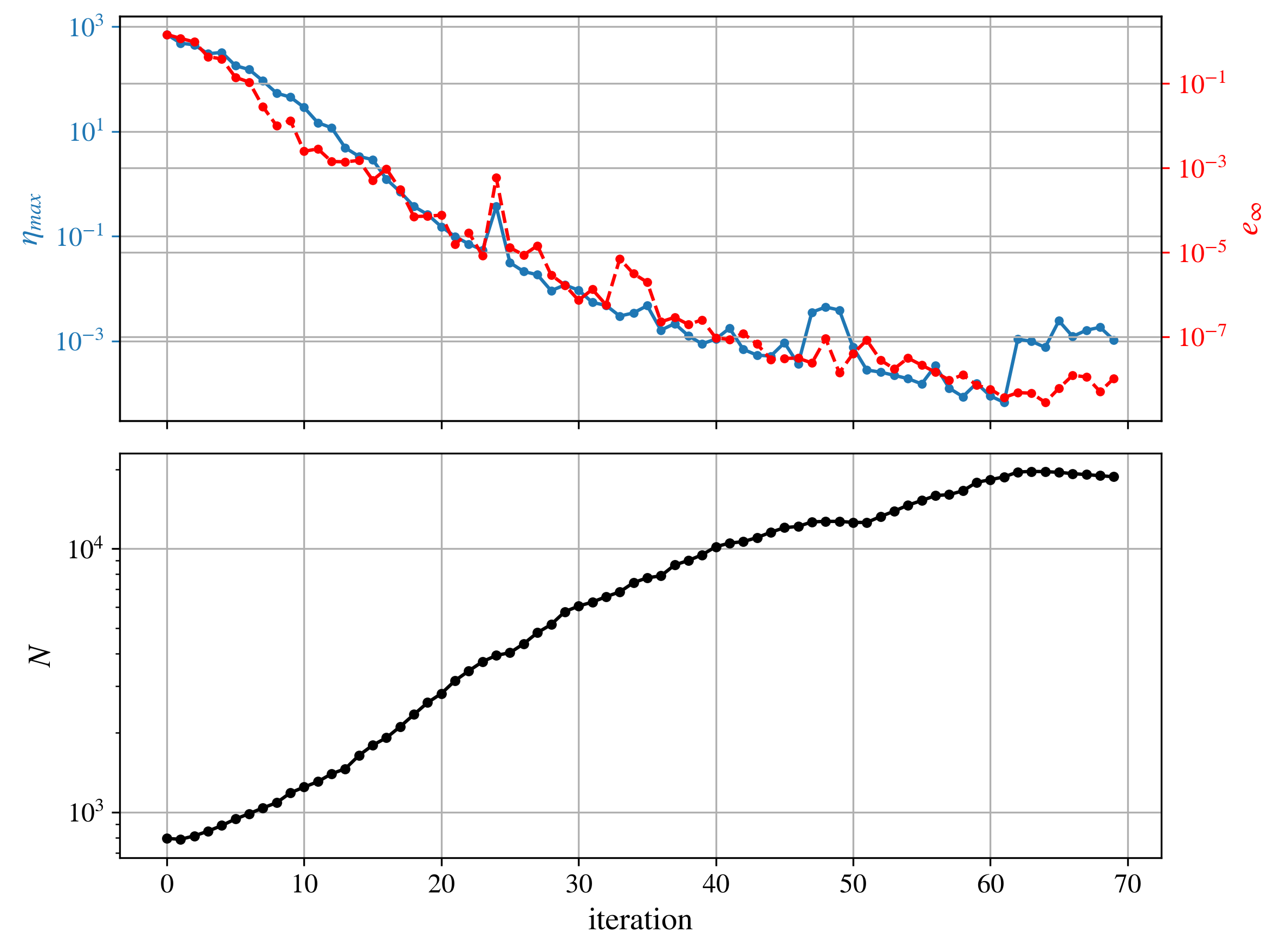}
    \caption{In the top row convergence of IMEX error indicator (blue) and convergence of numerical solution (red) within 70 iterations is shown, while the total number of computational nodes is shown below.}
    \label{fig:imex_convergence}
\end{figure}

The behaviour of IMEX over 70 adaptivity iterations is also studied in Figure~\ref{fig:imex_convergence}. We are pleased to find that the convergence limit of the indicator around iteration $N_{iter}=60$ agrees well with the convergence limit of the numerical solution.
This observation also makes the IMEX error indicator suitable for stopping criteria. Note that, in the process, the maximum error of the numerical solution has been reduced by about 9 orders of magnitude, while the maximum value of the error indicator has been reduced by about 7 orders of magnitude. In addition, Figure~\ref{fig:imex_convergence} also shows the number of computational nodes with respect to the adaptivity iterations.

\subsubsection{Approximation order distribution}
\label{sec:order_distribution}
The iterative adaptive procedure starts by obtaining the numerical solution of the unrefined problem setup. In this step, the approximation with the lowest approximation order, i.e.\ $m=2$, is assigned to all computational nodes. Later, the approximation orders are changed according to the marking and refinement strategies. Figure~\ref{fig:ref_imex_demonstration} shows the approximation order distributions for 3 selected adaptivity iterations. We can observe that the highest approximation orders are all near the exponentially strong source. Moreover, due to \h-adaptivity, the node density in the neighbourhood of the strong source is also significantly increased, i.e.\ $h_{max} / h_{min} \approx 52$ in the best-performing iteration.

After applying the \p-refinement strategy in the refinement step, the approximation order in two neighbouring nodes may differ by more than one. While numerical experiments with FEM have shown that heterogeneity of polynomial order in FEM leads to undesired oscillations of the approximated solution~\cite{wakeni2021p}, no similar behaviour was observed in our analyses with our setup using mesh-free methods. Thus, in contrast to \p-FEM, where additional smoothing of the approximation order takes place within the refinement module, we have completely avoided such manipulations and allow the approximation order in two neighbouring nodes to differ by more than one.



%
%
\subsubsection{Convergence rates of \hp-adaptive solution procedure}
\label{sec:conv_rates}

\begin{figure}
    \centering
    \includegraphics[width=0.9\textwidth]{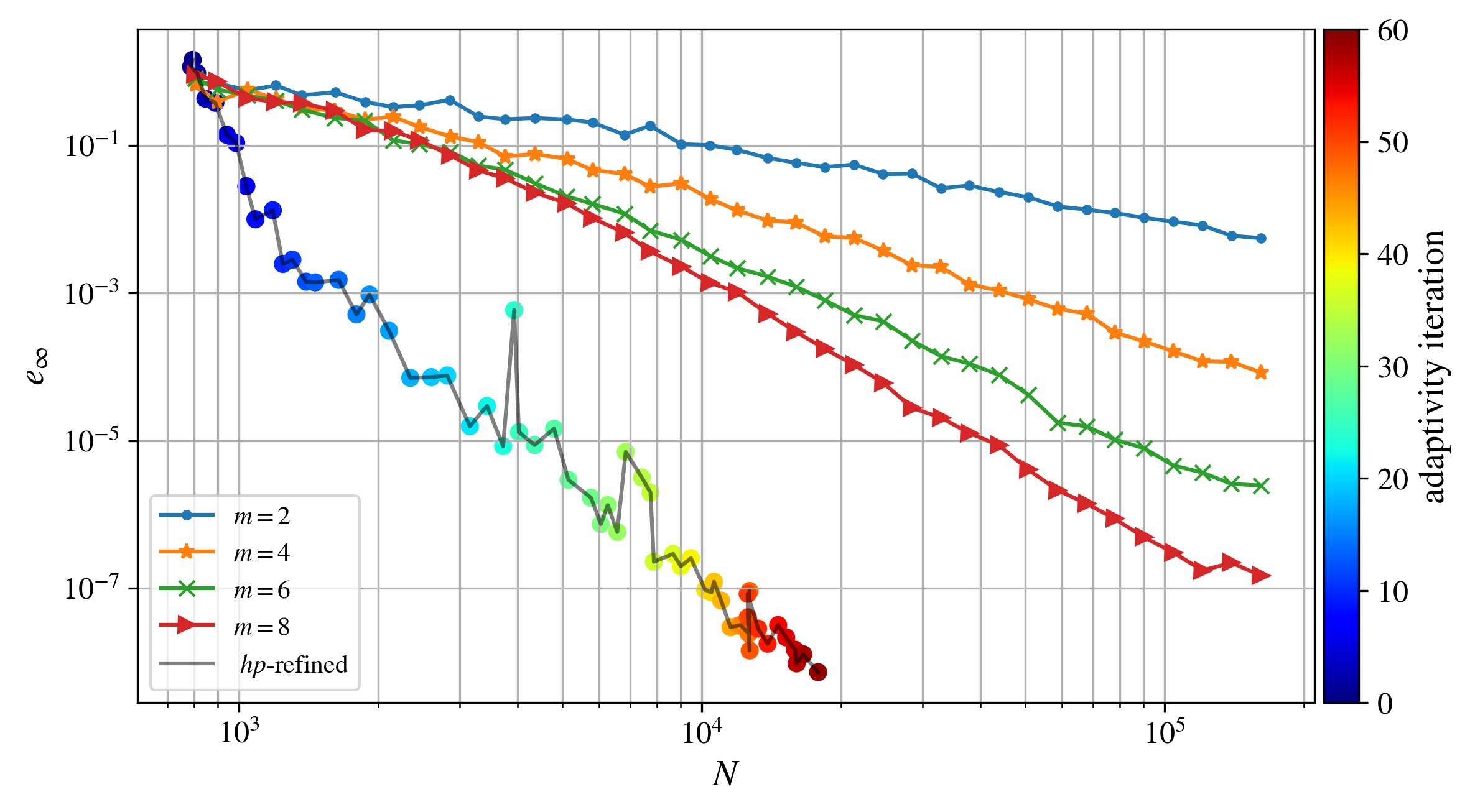}
    \caption{Convergence of the \hp-refined solution compared to the convergence of the unrefined solutions.}
    \label{fig:conv_2d}
\end{figure}

Finally, the convergence behaviour of the proposed \hp-adaptive solution procedure is studied. In addition to the convergence of a single \hp-adaptive run, Figure~\ref{fig:conv_2d} shows the convergences obtained without the use of refinement procedures, i.e.\ solutions obtained with uniform internodal spacing and approximation orders over the entire domain. The figure clearly shows that a \hp-adaptive solution procedure was able to significantly improve the numerical solution in terms of accuracy and computational points required.

As previously discussed by Eibner~\cite{eibner2007adaptive} and Demkowicz~\cite{demkowicz2002fully}, we believe that a more complex marking and refinement strategies would further improve the convergence behaviour, but already the proposed \hp-adaptive solution procedure significantly outperforms the unrefined solutions. Specifically, the refined solution is almost 4 orders of magnitude more accurate than the unrefined solution (for the highest approximation order $m=8$ used) at about $10^4$ computational nodes.





%
%

%
%
\section{Application to linear elasticity problems}
\label{sec:elasticity}
In this section we address two problems from linear elasticity that are conceptually different from the exponential peak problem discussed in Section~\ref{sec:demonstration}. While the solution of exponential peak problem is infinitely smooth, these two problems both have a singularity in the solution. 

In areas of smooth solution, the \hp-strategy should favour \p-refinement (assuming that the local discretization is sufficient, as briefly discussed in Section~\ref{sec:unrefined}), while near the singularity, \h-refinement should be preferred~\cite{mitchell2014comparison,demkowicz2002fully}.
However, the Texas Three Step based marking strategy used in this paper cannot trivially achieve this, since the strategy has no knowledge of the smoothness of the solution field. In addition, the strategy also cannot perform pure
\h- or pure \p-refinement~\cite{eibner2007adaptive} (see Figure~\ref{fig:indicator_complex}), which would be ideal in the limiting situations.
Instead, the strategy used enforces an increase in the approximation order by its design -- even if the solution is not smooth and even if low-regularity data is being used to construct the approximation. Nevertheless, in our experiments we observed an increase of the approximation order near the singularity only in the first few iterations, while the following iterations were focused on improving the local field description with \h-refinement. This observation is also in agreement with reports from the literature~\cite{mitchell2014comparison, eibner2007adaptive}, where authors justify the use of the Texas Three Step marking strategy also for problems with singularity in the solution.

\subsection{Fretting fatigue contact}
\label{sec:fatigue}
The application of the proposed \hp-adaptive solution procedure is further expanded to study a linear elasticity problem. Specifically, we obtain a \hp-refined solution to fretting fatigue contact problem~\cite{kosec_weak_2019} for which no closed form solution is known.

\begin{figure}
    \centering
    \includegraphics[width=0.8\textwidth]{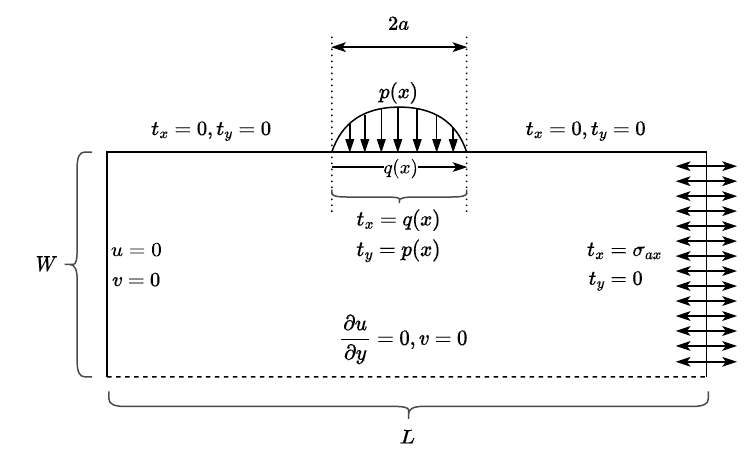}
    \caption{Fretting fatigue contact problem scheme and boundary conditions.}
    \label{fig:fwo_scheme}
\end{figure}

The problem dynamics is governed by the Cauchy-Navier equations
\begin{equation}
    (\lambda +\mu)\nabla (\nabla \cdot \b u)+\mu\nabla^2\b u = \b f
    \label{eq:NC}
\end{equation}
with unknown displacement vector $\b u$, external body force $\b f$ and Lam\'e parameters $\mu$ and $\lambda$. The domain of interest is a thin rectangle of width $W$, length $L$ and thickness $D$. Axial traction $\sigma_{ax}$ is applied to the right side of the rectangle, while a compression force is applied to the centre of the rectangle to simulate contact. The contact is simulated by a compressing force $F$ generated by two oscillating cylindrical pads of radius $R$, causing a tangential force $Q$. The tractions introduced by the two pads are predicted using an extension of Hertzian contact theory, which splits the contact area into the stick and slip zones depending on the friction coefficient $\mu$ and the combined elasticity modulus ${E^*}^{-1} = \Big ( \frac{1 - \nu_1^2}{E_1} + \frac{1 - \nu_2^2}{E_2} \Big)$, where $E_i$ and $\nu_i$ are the Young's modulus and the Poisson's ratios of the sample and the pad, respectively. The problem is shown schematically in Figure~\ref{fig:fwo_scheme} together with the boundary conditions.

Theoretical predictions from~\cite{slak2019adaptive} are used to obtain the contact half-width
\begin{align}
    a = 2 \sqrt{\frac{FR}{t\pi E^*}},
\end{align}
with normal traction
\begin{equation}
    p(x) = \begin{cases}
        p_0 \sqrt{1-\frac{x^2}{a^2}}, & |x| < a     \\
        0,                            & \text{else}
    \end{cases}, \qquad p_0 = \sqrt{\frac{FE^*}{t \pi R}},
\end{equation}
and tangential traction
\begin{equation}
    q(x) = \begin{cases}
        -\mu p_0 \left(\sqrt{1 - \frac{x^2}{a^2}} - \frac{c}{a}\sqrt{1 - \frac{(x-e)^2}{c^2}}\right), &
        |x-e| < c                                                                                                                                \\
        -\mu p_0 \sqrt{1 - \frac{x^2}{a^2}},                                                          & c \leq | x - e | \text{ and } |x| \leq a \\
        0,                                                                                            & \text{else}
    \end{cases}
\end{equation}
for $c = a\sqrt{1 - \frac{Q}{\mu f}}$ defined as the half-width of the slip zone and $e = \operatorname{sgn}(Q)\frac{a \sigma_{ax}}{4 \mu p_0}$ is the eccentricity due to axial loading. Note that the inequalities $Q \leq \mu F$ and $\sigma_{ax} \leq 4\left(1 - \sqrt{1 - \frac{Q}{\mu F}}\right)$ must hold for these expressions to be valid.

Plane strain approximation is used to reduce the problem from three to two dimensions and symmetry along the horizontal axis is used to further halve the problem size. Finally, $\Omega = [-L/2, L/2] \times [-W/2, 0]$ is taken as the domain.

We take $E_1 = E_2 = \unit[72.1]{GPa}$, $\nu_1 = \nu_2 = 0.33$, $L = \unit [40]{mm}$, $W = \unit [10]{mm}$, $t = \unit [4]{mm}$, $F = \unit [543]{N}$, $Q = \unit [155]{N}$, $\sigma_{ax} = \unit [100]{MPa}$, $R = \unit [10]{mm}$ and $\mu = 0.3$ for the model parameters. With this setup, the half-contact width $a$ is equal to $\unit[0.2067]{mm}$, which is about 200 times smaller than the domain width $W$. For stability reasons, the 4 corner nodes were removed after the domain was discretised.
\begin{table}
    \centering
    \renewcommand{\arraystretch}{1.2}
    \begin{tabular}{cccc|cccc||ccc} \hline
        \multicolumn{1}{c}{$\beta_h$} & \multicolumn{1}{c}{$\alpha_h$} & \multicolumn{1}{c}{$\lambda_h$} & \multicolumn{1}{c|}{$\vartheta_h$} & \multicolumn{1}{c}{$\beta_p$} & \multicolumn{1}{c}{$\alpha_p$} & \multicolumn{1}{c}{$\lambda_p$} & \multicolumn{1}{c||}{$\vartheta_p$} & \multicolumn{1}{c}{$h_{max}$} & \multicolumn{1}{c}{$N_{max}$} & \multicolumn{1}{c}{$N_{iter}$} \\ \hline
        $5\cdot 10^{-5}$              & $10^{-4}$                      & 5                               & 1.05                               & $10^{-3}$                     & 0.1                            & 4                               & 1.05                                & $2.5\cdot 10^{-4}$            & $5\cdot 10^{5}$               & 19                             \\ \hline
    \end{tabular}
    \caption{Adaptivity parameters used to obtain solution to fretting fatigue contact problem.}
    \label{tab:fwo_params}
\end{table}

The linear differential operators are approximated with RBF-FD using the PHS with order $k=3$ and local monomial augmentation limited to choose between approximation orders $m\in \left \{2, 4, 6, 8 \right \}$. The PHS order was left constant during the adaptive refinement. The \hp-refinement parameters used to obtain the numerical solution are given in Table~\ref{tab:fwo_params}.

\begin{figure}
    \centering
    \includegraphics[width=0.9\textwidth]{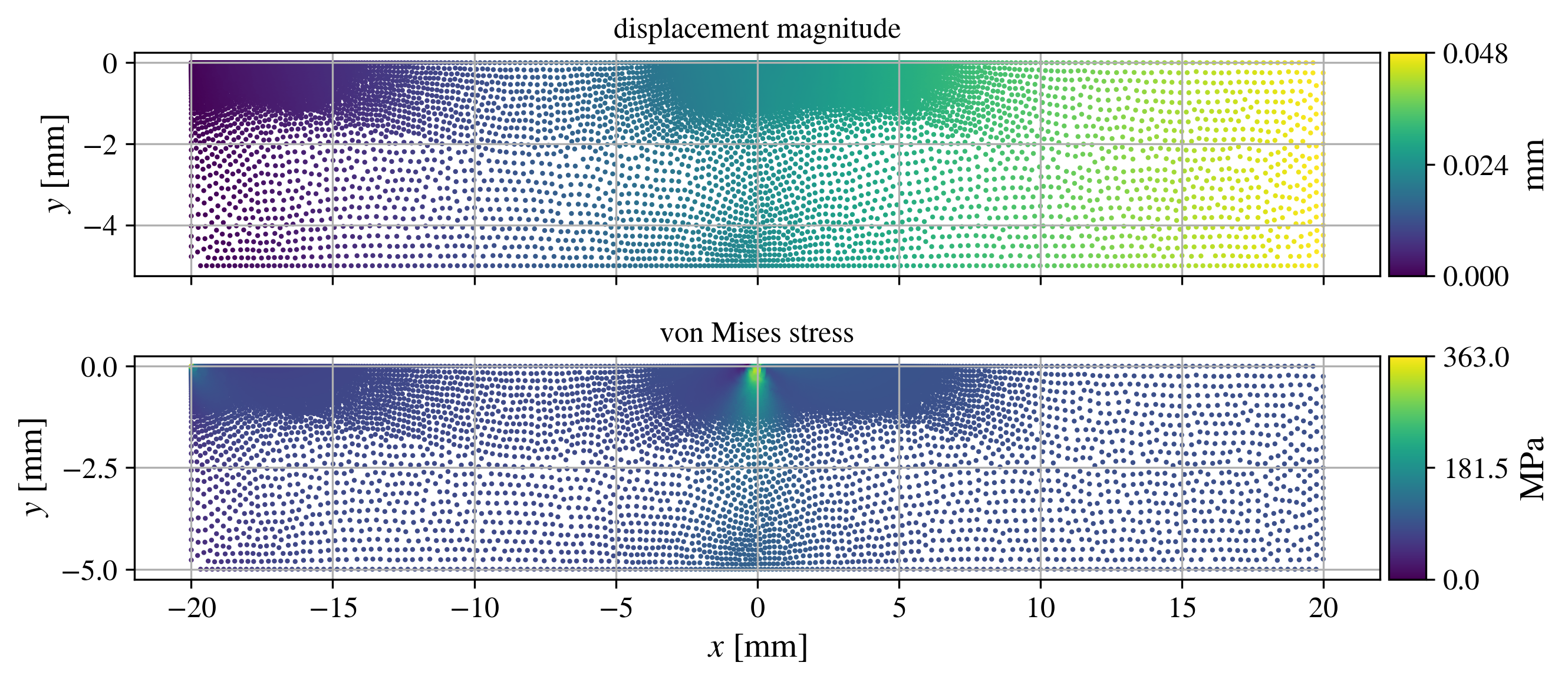}
    \caption{Example \hp-refined fretting fatigue contact solution.}
    \label{fig:fwo_example_solution}
\end{figure}

Figure~\ref{fig:fwo_example_solution} shows an example of a \hp-refined solution to fretting fatigue problem in the last adaptivity iteration with $N=46\,626$ computational nodes. We see that the solution procedure has successfully located the two critical points, i.e.\ the fixed upper left corner with a stress singularity and the area in the middle of the upper edge where contact is simulated. Note that the highest stress values (about $2$ times higher) were calculated in the singularity in the upper left corner, but these nodes are not shown as our focus is shifted towards the area under the contact.

\subsubsection{Surface traction under the contact}
\begin{figure}
    \centering
    \includegraphics[width=\textwidth]{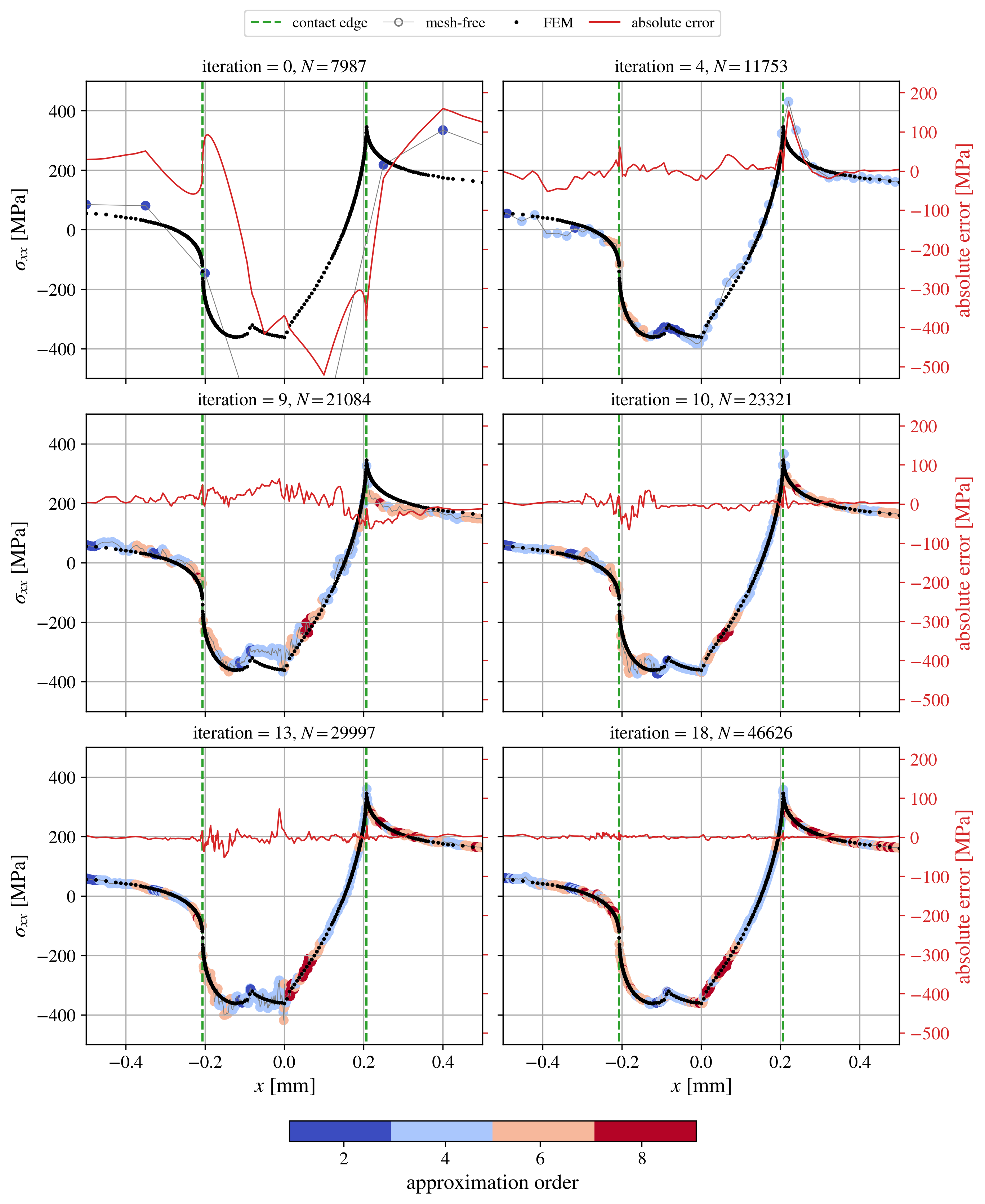}
    \caption{Surface traction under the contact for selected iteration steps demonstrating the \hp-adaptivity process. Colours are used to denote the local approximation orders. Numerical solution is additionally compared against the Abaqus FEM solution, where the red line is used to denote the absolute difference between the two methods. For clarity, the two dashed green lines show the edge contact.}
    \label{fig:fwo_conv}
\end{figure}

For a detailed analysis, we consider the surface traction $\sigma_{xx}$, as it is often used to determine the location of crack initiation. The surface traction is shown in Figure~\ref{fig:fwo_conv} for $6$ selected adaptivity iterations. The mesh-free nodes are coloured according to the local approximation order enforced by the \hp-adaptive solution procedure. The message of this figure is twofold. First, it is clear that the proposed IMEX error indicator can be successfully used in linear elasticity problems, and second, we find that the \hp-adaptive solution procedure has successfully approximated the surface traction near the contact. In doing so, the local field description under the contact has been significantly improved and the local approximation orders have taken a non-trivial distribution.

\begin{figure}
    \centering
    \includegraphics[width=0.8\textwidth]{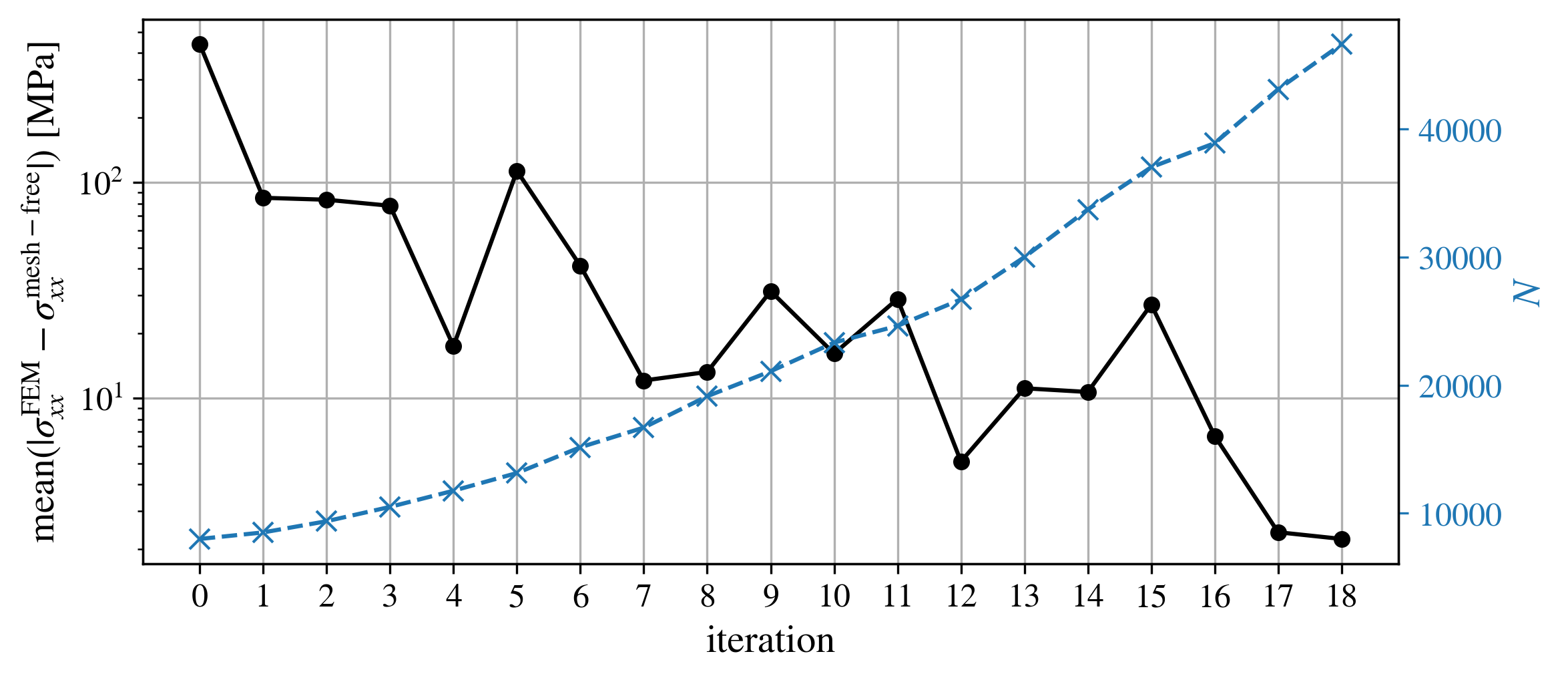}
    \caption{Mean surface traction difference between the two methods under the contact area.}
    \label{fig:fwo_conv_sxx}
\end{figure}

The surface traction in Figure~\ref{fig:fwo_conv} is additionally accompanied with the FEM results on a much denser mesh with more than 100\,000 DOFs obtained with the commercial solver Abaqus\textsuperscript{\textregistered}~\cite{kosec_weak_2019}. To calculate the absolute difference between the two methods, the mesh-free solution was interpolated to Abaqus's computational points using Sheppard's inverse distance weighting interpolation with $2$ nearest neighbours. We see that the absolute difference under the contact decreases with the number of adaptivity iterations and eventually settles at approximately 2 \% of the maximum difference from the initial iteration. As expected, the highest absolute difference is at the edges of the contact, i.e.\ around $x=a$ and $x=-a$, while the difference is even smaller in the rest of the contact area. The absolute difference between the two methods is further studied in Figure~\ref{fig:fwo_conv_sxx}, where the mean of $\lvert \sigma_{xx}^{ FEM } - \sigma_{xx}^{\text{mesh-free}}\rvert$ under the contact area, i.e.\ $-a \leq x \leq a$, is shown. We observe that the mesh-free \hp-refined solution converges towards the reference FEM solution with respect to the adaptivity iterations. Moreover, Figure~\ref{fig:fwo_conv_sxx} also shows the number of computational nodes with respect to the adaptivity iteration.

\subsection{The three-dimensional Boussinesq's problem}
\label{sec:boussinesq}
As a final benchmark problem we solve the three-dimensional Boussinesq's problem, where a concentrated normal traction acts on an isotropic half-space~\cite{Slaughter_2002}.

The problem has a closed form solution given in cylindrical coordinates $r$, $\theta$ and $z$ as
\begin{align}
    u_r         & = \frac{Pr}{4\pi \mu} \left(\frac{z}{R^3} - \frac{1-2\nu}{R(z+R)} \right), \qquad
    u_\theta = 0, \qquad
    u_z = \frac{P}{4\pi \mu} \left(\frac{2(1-\nu)}{R} + \frac{z^2}{R^3} \right), \nonumber          \\
    \sigma_{rr} & = \frac{P}{2\pi} \left( \frac{1-2\nu}{R(z+R)} - \frac{3r^2z}{R^5} \right), \qquad
    \sigma_{\theta\theta} = \frac{P(1-2\nu)}{2\pi} \left( \frac{z}{R^3} - \frac{1}{R(z+R)} \right),
    \label{eq:3d-problem}                                                                           \\
    \sigma_{zz} & = -\frac{3Pz^3}{2 \pi R^5}, \qquad
    \sigma_{rz} = -\frac{3Prz^2}{2 \pi R^5}, \qquad
    \sigma_{r\theta} = 0, \qquad \sigma_{\theta z} = 0, \nonumber
\end{align}
where $P$ is the magnitude of the concentrated force, $\nu$ is the Poisson's ratio, $\mu$ is the Lam\'e parameter and $R$ is the Eucledian distance to the origin. The solution has a singularity at the origin, which makes the problem ideal for treatment with adaptive procedures. Furthermore, the closed form solution also allows us to evaluate the accuracy of the numerical solution.

\begin{figure}
    \centering
    \includegraphics[width=0.4\textwidth]{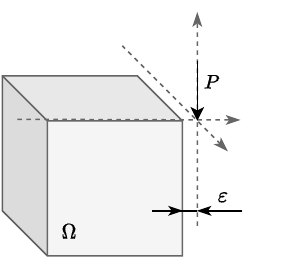}
    \caption{Schematic presentation of Boussinesq's problem.}
    \label{fig:bou_scheme}
\end{figure}

In our setup, we consider only a small part of the problem, i.e.\ $ \varepsilon = 0.1$ away from the singularity, as schematically shown in Figure~\ref{fig:bou_scheme}. From a numerical point of view, we solve the Navier-Cauchy Equation~\eqref{eq:NC} with Dirichlet boundary conditions described in~\eqref{eq:3d-problem}, where the domain $\Omega$ is defined as a box, i.e.\ $\Omega = [-1, -\varepsilon] \times [-1, -\varepsilon] \times [-1, -\varepsilon]$.

Although the closed form solution is given in cylindrical coordinate systems, the problem is implemented using cartesian coordinates. We employ the proposed mesh-free \hp-adaptive solution procedure where the shapes are computed with RBF-FD using the PHS with order $k=3$ and monomial augmentation restricted to choose between approximation orders $m\in \left \{2, 4, 6, 8 \right \}$. Other \hp-refinement related parameters are given in Table~\ref{tab:bou_params}. For the physical parameters of the problem, the values $P=-1$, $E=1$ and $\nu = 0.33$ were assumed. 

\begin{table}
    \centering
    \renewcommand{\arraystretch}{1.2}
    \begin{tabular}{cccc|cccc||ccc} \hline
        \multicolumn{1}{c}{$\beta_h$} & \multicolumn{1}{c}{$\alpha_h$} & \multicolumn{1}{c}{$\lambda_h$} & \multicolumn{1}{c|}{$\vartheta_h$} & \multicolumn{1}{c}{$\beta_p$} & \multicolumn{1}{c}{$\alpha_p$} & \multicolumn{1}{c}{$\lambda_p$} & \multicolumn{1}{c||}{$\vartheta_p$} & \multicolumn{1}{c}{$h_{max}$} & \multicolumn{1}{c}{$N_{max}$} & \multicolumn{1}{c}{$N_{iter}$} \\ \hline
        $10^{-3}$                     & $10^{-3}$                      & 3.75                            & 1.01                               & $10^{-4}$                     & $10^{-2}$                      & 3                               & 1.5                                 & 0.04                          & $7\cdot 10^{4}$               & 20                             \\ \hline
    \end{tabular}
    \caption{Adaptivity parameters used to obtain solution to Boussinesq's problem.}
    \label{tab:bou_params}
\end{table}

It is worth mentioning, that the final sparse system was solved using BiCGSTAB with ILUT preconditioner (employed with an initial guess obtained from the previous adaptivity iteration), where the global tolerance was set to $10^{-15}$ with a maximum number of $500$ iterations and drop-tolerance and fill-factor set to $10^{-6}$ and $60$ respectively. Other possible choices and their effect on the solution procedure are further discussed in Section~\ref{sec:sparse_system}.

\begin{figure}
    \centering
    \includegraphics[width=0.9\textwidth]{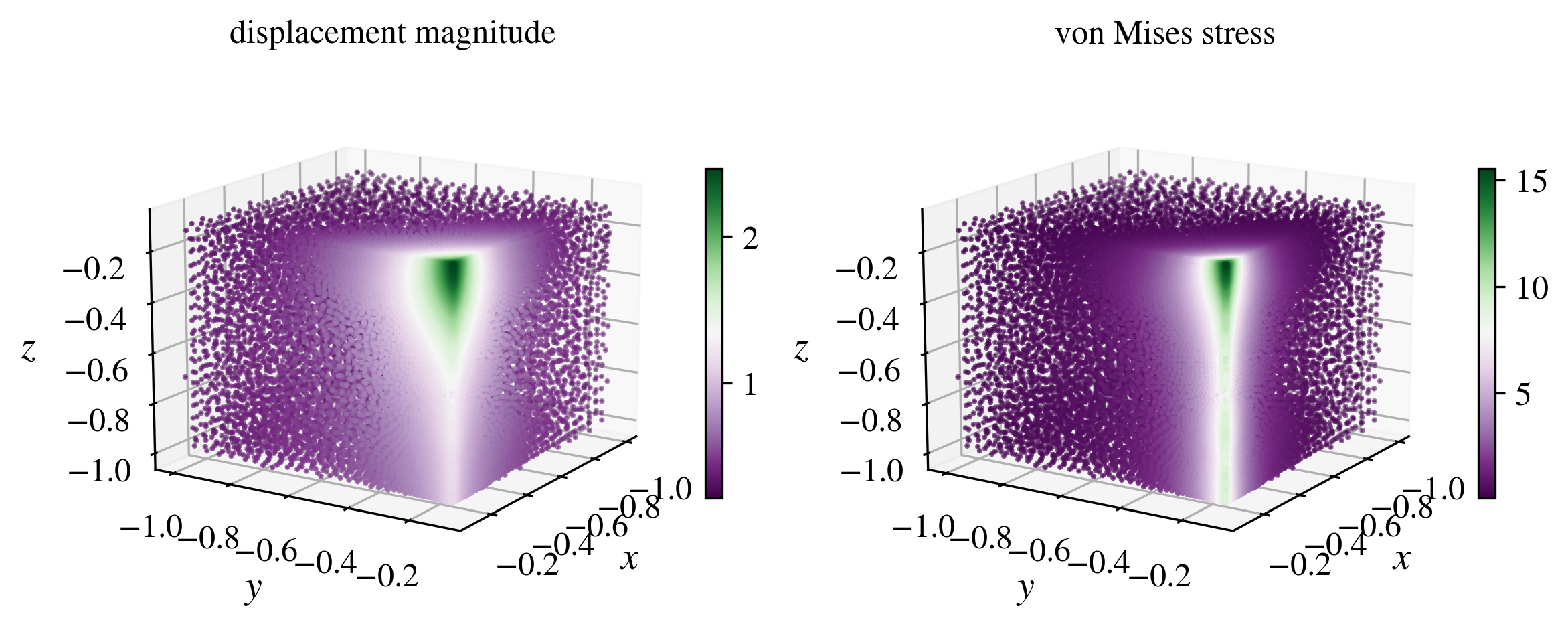}
    \caption{Example \hp-refined numerical solution to Boussinesq's problem.}
    \label{fig:bou_solution}
\end{figure}

Example \hp-refined numerical solution is given in Figure~\ref{fig:bou_solution}. We can see that the proposed \hp-adaptive solution procedure is sufficiently robust to obtain a good solution even for three-dimensional problems with singularities. Additionally, we also observe that the IMEX error indicator successfully identified the singularity, effectively seen as an increase in the local field description in the neighbourhood of the concentrated force applied at the origin.

\subsubsection{The von Mises stress along the body diagonal}
\begin{figure}
    \centering
    \includegraphics[width=0.9\textwidth]{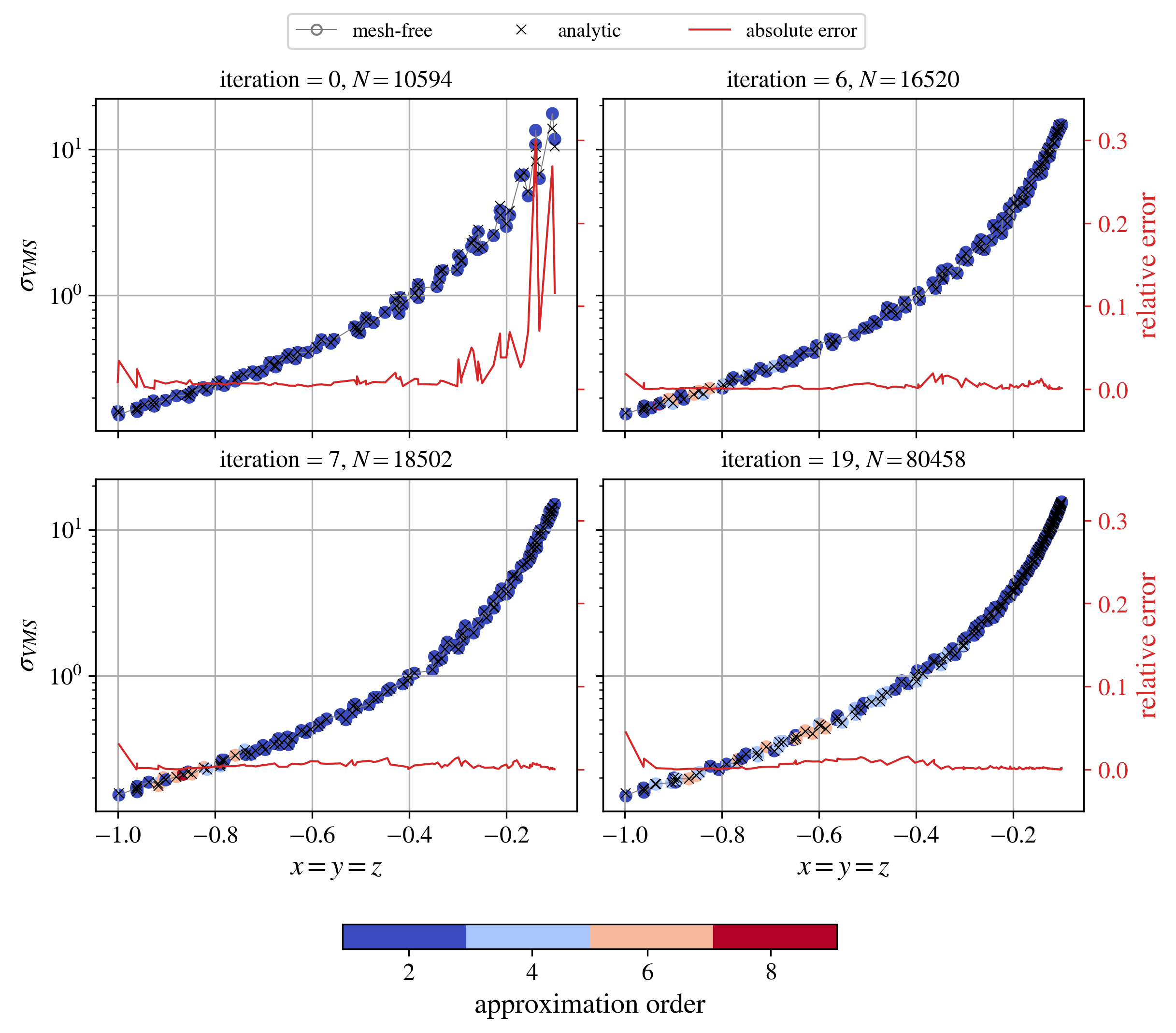}
    \caption{Numerical solution compared to analytical solution at the nodes near the body diagonal $(-1, -1, -1) \rightarrow (-\varepsilon, -\varepsilon, -\varepsilon)$ for selected iterations.}
    \label{fig:bou_iterations}
\end{figure}

Figure~\ref{fig:bou_iterations} shows further evaluation of the \hp-refined mesh-free numerical solution. Here, the von Mises stress at points near the body diagonal $(-1, -1, -1) \rightarrow (-\varepsilon, -\varepsilon, -\varepsilon)$ is calculated for selected 4 adaptivity iterations and compared to the analytical values in terms of relative error. In addition, the nodes are coloured according to the local approximation order enforced by the \hp-adaptive solution procedure. We can see that the highest relative error of approximately $0.3$ at the initial state is observed in the neighbourhood of the origin. In the final iteration, the relative error is reduced by about an order of magnitude. We also see that the \hp-adaptive solution procedure has found a non-trivial order distribution and that the number of nodes in the neighbourhood of the corner $(-\varepsilon, -\varepsilon, -\varepsilon)$ has increased significantly.

\begin{figure}
    \centering
    \includegraphics[width=0.9\textwidth]{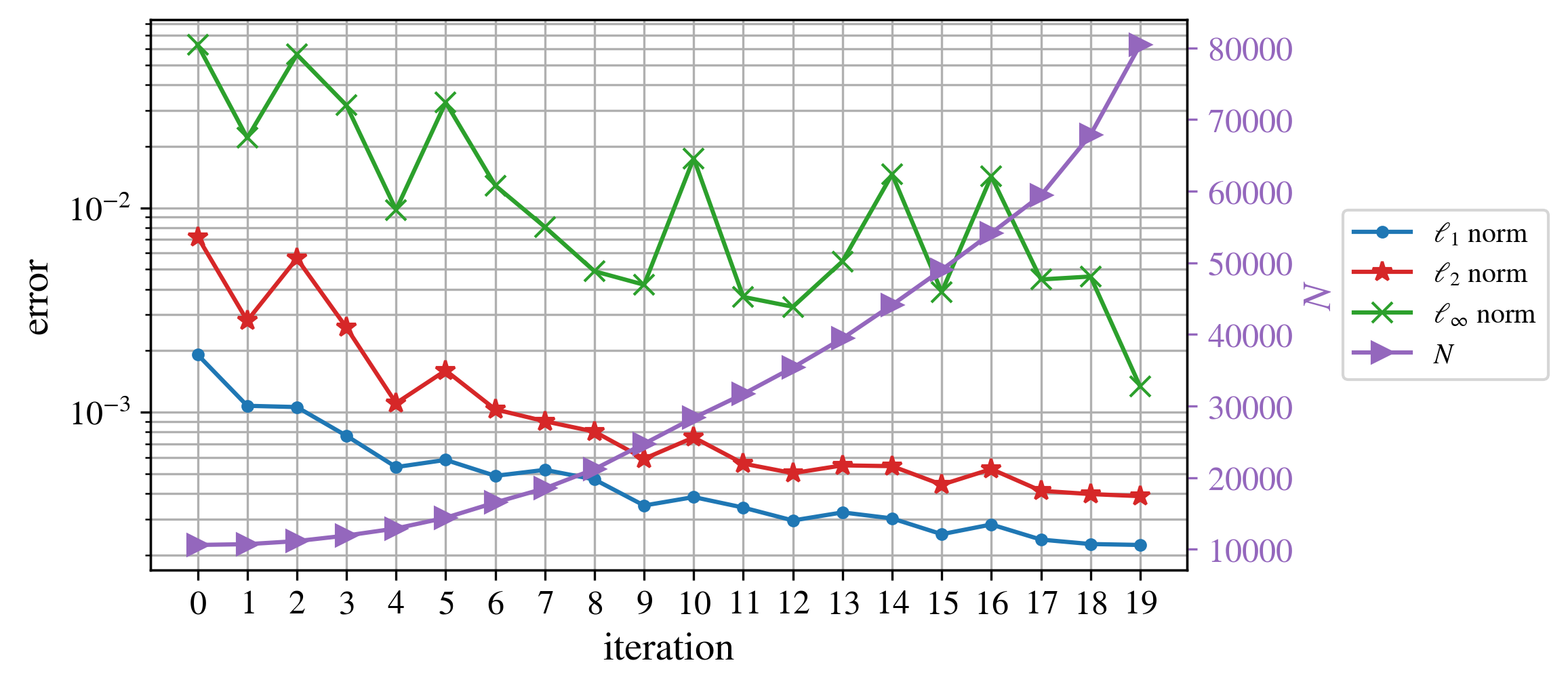}
    \caption{Convergence of numerical solution along with number of computational nodes.}
    \label{fig:conv_bou}
\end{figure}

A more quantitative analysis of the mesh-free solution is given in Figure~\ref{fig:conv_bou} where the $\ell _1$, $\ell_2$ and $\ell_{\infty}$ error norms and number of computational nodes vs.\ adaptivity iteration are shown. Compared to the initial state, the \hp-adaptive solution procedure was able to achieve a numerical solution that was almost two orders of magnitude more accurate. In the process, the number of computational nodes increased from 10\,500 in the initial state to about 80\,000 in the final iteration. However, it is interesting to observe that with the configuration from Table~\ref{tab:bou_params}, none of the computational nodes used the approximation with the highest order allowed ($m=8$). Instead, in the final iteration, there were $130$ nodes approximated with $m=6$, and $5937$ with $m=4$, while the rest were approximated with the second order. Note that, as expected, most of the higher order approximations are near the concentrated force -- which is difficult to represent visually, so we only give the descriptive data.

For reference, we take the \h-refined solution by Slak et~al.~\cite{slak2019adaptive}, who were able to reduce the infinity norm error by about an order of magnitude with $N\approx 140\,000$ nodes in the final iteration. It is perhaps naive to compare this result with ours, since the authors use different marking and refinement strategies and, more importantly, a different error indicator. Nevertheless, the infinity norm error of our \hp-refined solution is in the neighbourhood of $10^{-3}$ compared to theirs at approximately $10^{-2}$ with almost twice as many computational nodes. We believe our results could be further improved by fine-tuning the free parameters, but we decided to avoid such an approach.

\subsubsection{Additional discussion on solving the global sparse system}
\label{sec:sparse_system}
In all previous sections, we have completely neglected the importance of solving the global sparse system in the proposed \hp-adaptive solution procedure with a suitable solver. However, inappropriate choice of solver can lead to inaccurate or even unstable behaviour and, most importantly, unreasonably large computational cost. To avoid such flaws, we compared an iterative BiCGSTAB and BiCGSTAB with ILUT preconditioner with two direct solvers --- namely the SparseLU and the PardisoLU --- on a \hp-adaptive solution to the Boussinesq problem, performing $25$ adaptivity iterations with approximately $10\,000$ initial nodes and $135\,000$ nodes after the last iteration. Note that the iterative BiCGSTAB solver with ILUT preconditioner was employed with an initial guess obtained from the previous adaptivity iteration.

In addition to the discussed solvers, we also tried the SparseQR. While its stability and accuracy were comparable to other solvers, its computational cost was significantly higher and was therefore removed from further analysis and from the list of potential candidates. For all performed tests we used the EIGEN linear algebra library~\cite{eigenweb}.

\begin{figure}
    \centering
    \includegraphics[width=0.98\textwidth]{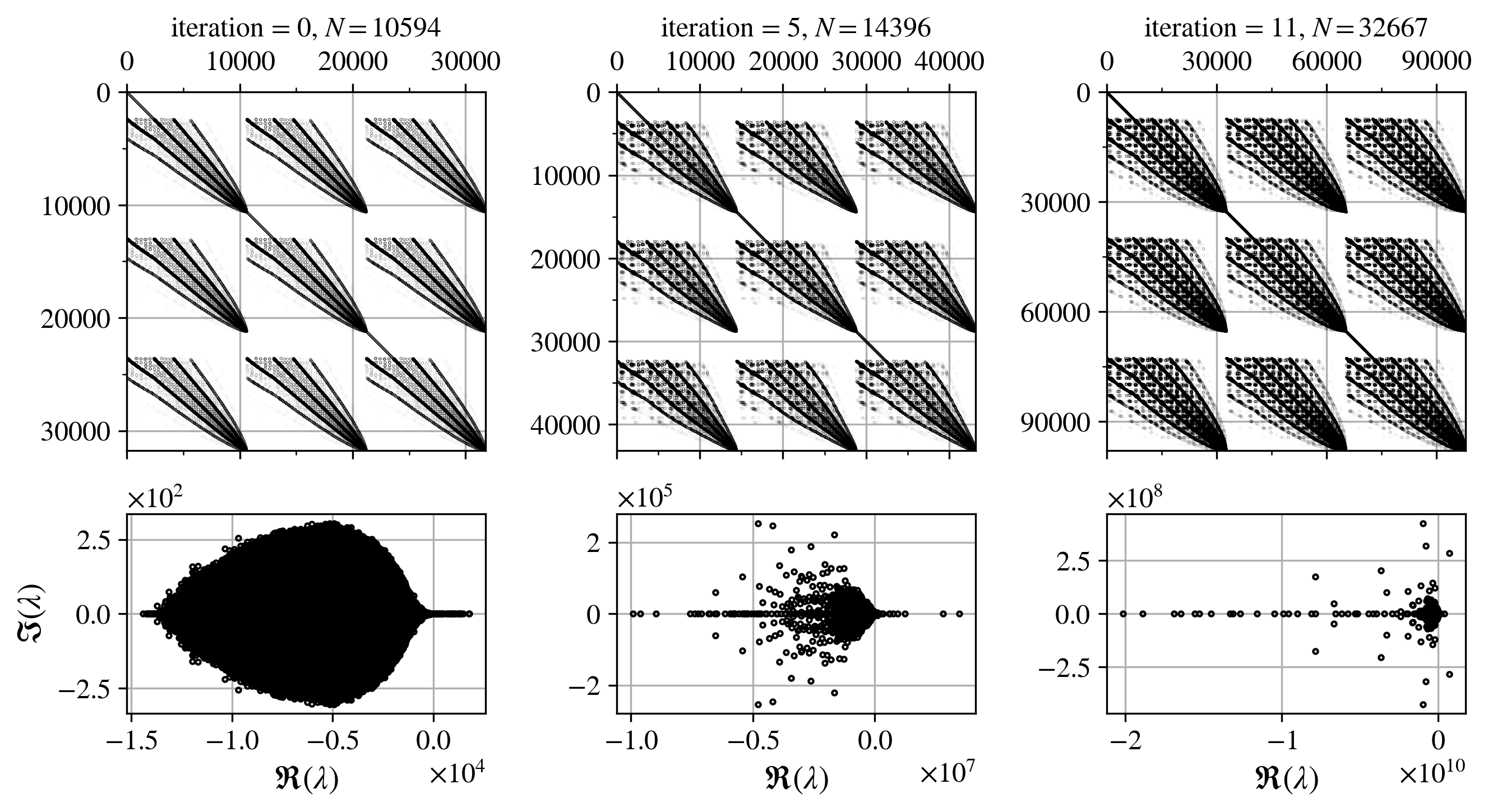}
    \caption{Global sparse matrix plot (top row) and spectra of the matrices (bottom row) at three selected iterations of the \hp-adaptive solution procedure. Note that the spectra are computed for the BiCGSTAB solver with an ILUT preconditioning using an estimate from the previous iteration.}
    \label{fig:spectra}
\end{figure}

Let us first examine the sparsity patterns of the systems assembled at different stages of the \hp-adaptive process in Figure~\ref{fig:spectra}, where we can see how the system increases in size and also becomes less sparse due to globally decreasing the internodal distance \h~and increasing the approximation order \p. Additionally, the spectra of the matrices are shown in the bottom row of Figure~\ref{fig:spectra}, where we can see that the ratios between the real and imaginary parts of the eigenvalues are in good agreement with previous studies~\cite{janvcivc2021monomial,bayona2017role,bayona2019comparison}.

\begin{figure}
    \centering
    \includegraphics[width=\textwidth]{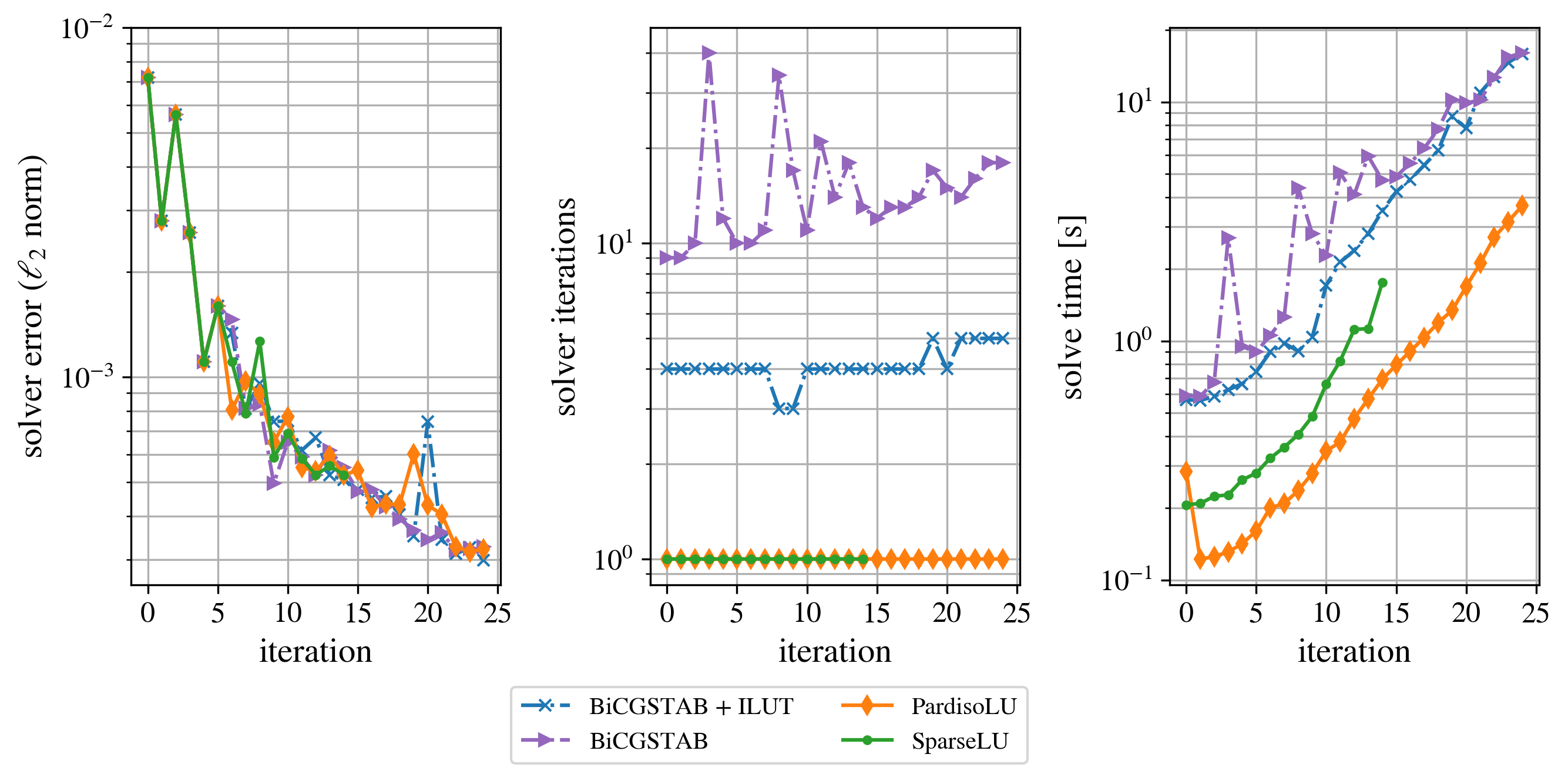}
    \caption{Error of the final solution with respect to the adaptivity iteration for different solvers (left), number of solver iteration per adaptivity iteration (centre) and solver compute time for each adaptivity iteration (right).}
    \label{fig:solvers}
\end{figure}

Moreover, Figure~\ref{fig:solvers} presents three different views of the solvers' performance: (i) the achieved accuracy of the final solution for different solvers, (ii) the number of iterations a solver needs to converge, and (iii) the execution times of each solver, each with respect to the \hp-adaptive iterations. The differences in final accuracy for different solvers are marginal. Perhaps the BiCGSTAB shows better stability behaviour (in terms of error scatter) compared to others. Nevertheless, it is important to observe, that the SparseLU only works until the $15$th iteration with approximately 50\,000 nodes, at which point our computer ran out of the available 12 Gb memory, which is to be expected due to the computational complexity or SparseLU. PardisoLU, on the other hand, remains stable through all adaptivity iterations.

Generally speaking, the number of iterations BiCGSTAB needs to converge increases with \hp-adaptivity iterations due to the increasing non-zero elements in the global system. The BiCGSTAB with a ILUT precoditioner shows similar behaviour, but with approximately $2/3$ less iterations required. Both direct solvers, of course, require only one ``iteration''. Finally, the analysis of the execution time shows that the PardisoLU solver is by far the most efficient among all considered candidates.

With all things considered, PardisoLU seems to be the the best candidate for \hp-adaptive solution procedure. However, the last adaptivity iteration with approximately 115\,000 nodes was coincidentally right at the limit of the available 12 Gb RAM memory -- using approximately 10.5 Gb. It is therefore expected that like SparseLU, the PardisoLU would soon run out of memory for larger domains. To avoid such problems, we chose to work with a general purpose iterative BiCGSTAB solver with ILUT preconditioner employed with an initial guess, since it shows slightly better computational efficiency than the pure BiCGSTAB and required only 7.5 Gb of RAM for approximately 135\,000 nodes in the final adaptivity iteration.

\section{Conclusions}
\label{sec:conclusions}
In this paper we establish a baseline for \hp~strong form mesh-free analysis. We have formulated and implemented a \hp-adaptive solution procedure and demonstrated its performance in three different numerical experiments.

The cornerstone of the presented \hp-adaptive method is an iterative solve--estimate--mark--refine paradigm with the modified Texas Three Step marking strategy. The \h-refine of the proposed method relies on an advancing front node positioning algorithm based on Poisson disc sampling, which enables dimension-independent node generation supporting spatially variable density distributions. For the adaptive order of the method, we exploit an elegant control over the order of the approximation via the augmenting monomials in the approximation basis.

We proposed an IMEX error indicator, where the implicit solution of the problem is processed with the higher order local explicit representation of PDE at hand, e.g.\ if the implicit solution is computed with a second order approximation, the explicit re-evaluation happens at fourth order. Our analyses show that the proposed error indicator successfully captures main characteristics of error distributions, which suffices for the proposed iterative adaptivity.

The proposed \hp-adaptive solution procedure is first demonstrated on a two-dimensional Poisson problem with exponential source and mixed boundary conditions. Further demonstration focuses on linear elasticity problems. First, a 2D fretting fatigue problem -- a contact problem with pronounced peaks in the surface stress, and second, a 3D Boussinesq's problem with stress singularity. In both cases, we have demonstrated the advantages of using the proposed \hp-adaptive approach.

Although the \hp-adaptivity introduces additional steps in the solution procedure 
and is therefore undoubtedly computationally more expensive per node than the non-adaptive approach, it is essential in problems that exhibit volatilities in solution in small regions of the domain. For example singularity in the contact problem require excessively detailed numerical analysis near the contact compared to the rest (the bulk) of the domain. Such cases are extremely difficult (or even impossible) to solve without adaptivity, since the minimal uniform \h~and \p~distribution required to capture these volatilities would lead to unreasonably high computational complexity. In cases of a smooth solution, however, the benefits of \hp-adaptivity in most cases do not justify its computational overheads.

We are aware that there are many opportunities for improvement of presented methodology. The IMEX error indicator needs further clarification. Other error indicators should also be implemented and tested. During our experiments, we have found that a marking strategy with more free parameters leads to better accuracy, but is also more difficult to understand and control and can be case dependent. A smarter and more effective refinement and marking strategies are certainly part of future work. These should possibly take into account more information about the method itself, e.g. the dependence of the computational complexity on the approximation order, and most importantly local data regularity to choose between \p~and \h~refinement.

One of our goals in future work is also generalisation of the presented \hp-adaptive solution procedure to time-dependent problems. The most straightforward approach to achieve that is to granularly adapt \h~and \p~throughout the simulation. In its simplest form, the proposed \hp-adaptivity would be performed at each time step, starting with the \hp~distributions of the previous time step and using the same adaptivity parameters for all time steps. A more sophisticated approach would also take into account the desired accuracy during the simulation, resulting in time-dependent adaptivity parameters. For example, if one is only interested in a steady state solution, the desired accuracy would increase with time, reaching its maximum at steady state. Additionally, to perform proper adaptive analysis, the time step should also be adaptive, which requires an additional step in the \hp-adaptive solution procedure.


%
%
\section*{Acknowledgments}
The authors acknowledge the financial support from the Slovenian Research Agency research core funding No.\ P2-0095, and research projects No.\ J2-3048 and No.\ N2-0275 also funded by National Science Centre, Poland under the OPUS call in the Weave programme 2021/43/I/ST3/00228.

\section*{Declarations}
\textbf{Conflict of interest.} The authors declare that they have no conflict of interest. All the co-authors have confirmed to know the submission of the manuscript by the corresponding author.


%
%

\bibliography{bib}

\begin{thebibliography}{10}
\providecommand{\url}[1]{{#1}}
\providecommand{\urlprefix}{URL }
\providecommand{\doi}[1]{\url{https://doi.org/#1}}
\bibcommenthead

\bibitem{upadhyay2021numerical}
B.D. Upadhyay, S.S. Sonigra, S.D. Daxini, Numerical analysis perspective in
  structural shape optimization: A review post 2000.
\newblock Advances in Engineering Software \textbf{155}, 102,992 (2021)

\bibitem{mitchell2014comparison}
W.F. Mitchell, M.A. McClain, A comparison of hp-adaptive strategies for
  elliptic partial differential equations.
\newblock ACM Transactions on Mathematical Software (TOMS) \textbf{41}(1),
  1--39 (2014)

\bibitem{SEGETH20101589}
K.~Segeth, A review of some a posteriori error estimates for adaptive finite
  element methods.
\newblock Mathematics and Computers in Simulation \textbf{80}(8), 1589--1600
  (2010).
\newblock \doi{https://doi.org/10.1016/j.matcom.2008.12.019}.
\newblock
  \urlprefix\url{https://www.sciencedirect.com/science/article/pii/S0378475408004230}.
\newblock ESCO 2008 Conference

\bibitem{liu2005introduction}
G.R. Liu, Y.T. Gu, \emph{An introduction to meshfree methods and their
  programming} (Springer Science \& Business Media, 2005)

\bibitem{liu2002mesh}
G.R. Liu, \emph{Mesh free methods: moving beyond the finite element method}
  (CRC press, 2002).
\newblock \doi{10.1201/9781420040586}

\bibitem{zienkiewicz2005finite}
O.C. Zienkiewicz, R.L. Taylor, J.Z. Zhu, \emph{The finite element method: its
  basis and fundamentals} (Elsevier, 2005)

\bibitem{van2021fast}
K.~van~der Sande, B.~Fornberg, Fast variable density 3-d node generation.
\newblock SIAM Journal on Scientific Computing \textbf{43}(1), A242--A257
  (2021)

\bibitem{shankar2018robust}
V.~Shankar, R.M. Kirby, A.L. Fogelson, Robust node generation for mesh-free
  discretizations on irregular domains and surfaces.
\newblock SIAM Journal on Scientific Computing \textbf{40}(4), A2584--A2608
  (2018)

\bibitem{jacquemin_smart_2023}
T.~Jacquemin, P.~Suchde, S.P. Bordas, Smart {Cloud} {Collocation}:
  {Geometry}-{Aware} {Adaptivity} {Directly} {From} {CAD}.
\newblock Computer-Aided Design \textbf{154}, 103,409 (2023).
\newblock \doi{10.1016/j.cad.2022.103409}.
\newblock
  \urlprefix\url{https://linkinghub.elsevier.com/retrieve/pii/S0010448522001427}

\bibitem{slak2019adaptive}
J.~Slak, G.~Kosec, Adaptive radial basis function--generated finite differences
  method for contact problems.
\newblock International Journal for Numerical Methods in Engineering
  \textbf{119}(7), 661--686 (2019)

\bibitem{davydov2011adaptive}
O.~Davydov, D.T. Oanh, Adaptive meshless centres and rbf stencils for poisson
  equation.
\newblock Journal of Computational Physics \textbf{230}(2), 287--304 (2011)

\bibitem{jacquemin_unified_2021}
T.~Jacquemin, S.P.A. Bordas, A unified algorithm for the selection of
  collocation stencils for convex, concave, and singular problems.
\newblock International Journal for Numerical Methods in Engineering
  \textbf{122}(16), 4292--4312 (2021)

\bibitem{janvcivc2021monomial}
M.~Jan{\v{c}}i{\v{c}}, J.~Slak, G.~Kosec, Monomial augmentation guidelines for
  rbf-fd from accuracy versus computational time perspective.
\newblock Journal of Scientific Computing \textbf{87}(1), 1--18 (2021)

\bibitem{bayona2017role}
V.~Bayona, N.~Flyer, B.~Fornberg, G.A. Barnett, On the role of polynomials in
  rbf-fd approximations: Ii. numerical solution of elliptic pdes.
\newblock Journal of Computational Physics \textbf{332}, 257--273 (2017)

\bibitem{belytschko1996meshless}
T.~Belytschko, Y.~Krongauz, D.~Organ, M.~Fleming, P.~Krysl, Meshless methods:
  an overview and recent developments.
\newblock Computer methods in applied mechanics and engineering
  \textbf{139}(1-4), 3--47 (1996)

\bibitem{kosec_simulation_2014}
G.~Kosec, B.~Šarler, Simulation of macrosegregation with mesosegregates in
  binary metallic casts by a meshless method.
\newblock Engineering Analysis with Boundary Elements \textbf{45}, 36--44
  (2014).
\newblock \doi{10.1016/j.enganabound.2014.01.016}.
\newblock
  \urlprefix\url{https://linkinghub.elsevier.com/retrieve/pii/S0955799714000290}

\bibitem{maksic_cooling_2019}
M.~Maksić, V.~Djurica, A.~Souvent, J.~Slak, M.~Depolli, G.~Kosec, Cooling of
  overhead power lines due to the natural convection.
\newblock International Journal of Electrical Power \& Energy Systems
  \textbf{113}, 333--343 (2019).
\newblock \doi{10.1016/j.ijepes.2019.05.005}.
\newblock
  \urlprefix\url{https://linkinghub.elsevier.com/retrieve/pii/S0142061518340055}

\bibitem{gui1985h}
W.z. Gui, I.~Babuska, The h, p and hp versions of the finite element method in
  1 dimension. part 3. the adaptive hp version.
\newblock Tech. rep., MARYLAND UNIV COLLEGE PARK LAB FOR NUMERICAL ANALYSIS
  (1985)

\bibitem{gui1986h}
W.Z. Gui, I.~Babu{\v{s}}ka, The h, p and hp versions of the finite element
  method in 1 dimension. part ii. the error analysis of the h and hp versions.
\newblock Numerische Mathematik \textbf{49}(6), 613--657 (1986)

\bibitem{devloo12recent}
P.R. Devloo12, C.M. Bravo, E.C. Rylo, Recent developments in hp adaptive
  refinement

\bibitem{zienkiewicz1987simple}
O.C. Zienkiewicz, J.Z. Zhu, A simple error estimator and adaptive procedure for
  practical engineerng analysis.
\newblock International journal for numerical methods in engineering
  \textbf{24}(2), 337--357 (1987)

\bibitem{gonzalez2021error}
O.A. Gonz{\'a}lez-Estrada, S.~Natarajan, J.J. R{\'o}denas, S.P. Bordas, Error
  estimation for the polygonal finite element method for smooth and singular
  linear elasticity.
\newblock Computers \& Mathematics with Applications \textbf{92}, 109--119
  (2021)

\bibitem{Ebrahimnejad}
M.~Thimnejad, N.~Fallah, A.R. Khoei, Adaptive refinement in the meshless finite
  volume method for elasticity problems.
\newblock Computers \& Mathematics with Applications \textbf{69}(12),
  1420--1443 (2015).
\newblock \doi{10.1016/j.camwa.2015.03.023}

\bibitem{Angulo}
A.~Angulo, L.P. Pozo, F.~Perazzo, A posteriori error estimator and an adaptive
  technique in meshless finite points method.
\newblock Engineering Analysis with Boundary Elements \textbf{33}(11),
  1322--1338 (2009).
\newblock \doi{10.1016/j.enganabound.2009.06.004}

\bibitem{oanh2017adaptive}
D.T. Oanh, O.~Davydov, H.X. Phu, Adaptive rbf-fd method for elliptic problems
  with point singularities in 2d.
\newblock Applied Mathematics and Computation \textbf{313}, 474--497 (2017)

\bibitem{Sang}
P.~Sang‐Hoon, K.~Kie‐Chan, Y.~Sung-Kie, A posterior error estimates and an
  adaptive scheme of least-squares meshfree method.
\newblock International Journal for Numerical Methods in Engineering
  \textbf{58}(8), 1213--1250 (2003).
\newblock \doi{10.1002/nme.817}

\bibitem{afshar2011node}
M.~Afshar, M.~Naisipour, J.~Amani, Node moving adaptive refinement strategy for
  planar elasticity problems using discrete least squares meshless method.
\newblock Finite Elements in Analysis and Design \textbf{47}(12), 1315--1325
  (2011)

\bibitem{guo1986hp}
B.~Guo, I.~Babu{\v{s}}ka, The hp version of the finite element method.
\newblock Computational Mechanics \textbf{1}(1), 21--41 (1986)

\bibitem{mitchell2016performance}
W.F. Mitchell, Performance of hp-adaptive strategies for 3d elliptic problems
  (2016)

\bibitem{tinsley1995three}
J.~Tinsley~Oden, W.~Wu, M.~Ainsworth, in \emph{Modeling, mesh generation, and
  adaptive numerical methods for partial differential equations} (Springer,
  1995), pp. 347--366

\bibitem{ainsworth1997aspects}
M.~Ainsworth, B.~Senior, Aspects of an adaptive hp-finite element method:
  Adaptive strategy, conforming approximation and efficient solvers.
\newblock Computer Methods in Applied Mechanics and Engineering
  \textbf{150}(1-4), 65--87 (1997)

\bibitem{houston2005note}
P.~Houston, B.~Senior, E.~S{\"u}li, in \emph{Numerical mathematics and advanced
  applications} (Springer, 2003), pp. 631--656

\bibitem{houston2005note1}
P.~Houston, E.~S{\"u}li, A note on the design of hp-adaptive finite element
  methods for elliptic partial differential equations.
\newblock Computer Methods in Applied Mechanics and Engineering
  \textbf{194}(2-5), 229--243 (2005)

\bibitem{eibner2007adaptive}
T.~Eibner, J.M. Melenk, An adaptive strategy for hp-fem based on testing for
  analyticity.
\newblock Computational Mechanics \textbf{39}(5), 575--595 (2007)

\bibitem{burg2011convergence}
M.~B{\"u}rg, W.~D{\"o}rfler, Convergence of an adaptive hp finite element
  strategy in higher space-dimensions.
\newblock Applied numerical mathematics \textbf{61}(11), 1132--1146 (2011)

\bibitem{demkowicz2002fully}
L.~Demkowicz, W.~Rachowicz, P.~Devloo, A fully automatic hp-adaptivity.
\newblock Journal of Scientific Computing \textbf{17}(1), 117--142 (2002)

\bibitem{rachowicz2006fully}
W.~Rachowicz, D.~Pardo, L.~Demkowicz, Fully automatic hp-adaptivity in three
  dimensions.
\newblock Computer methods in applied mechanics and engineering
  \textbf{195}(37-40), 4816--4842 (2006)

\bibitem{benito2003h}
J.~Benito, F.~Urena, L.~Gavete, R.~Alvarez, An h-adaptive method in the
  generalized finite differences.
\newblock Computer methods in applied mechanics and engineering
  \textbf{192}(5-6), 735--759 (2003)

\bibitem{liu2006stabilized}
G.~Liu, B.B. Kee, L.~Chun, A stabilized least-squares radial point collocation
  method (ls-rpcm) for adaptive analysis.
\newblock Computer methods in applied mechanics and engineering
  \textbf{195}(37-40), 4843--4861 (2006)

\bibitem{hu2019spatially}
W.~Hu, N.~Trask, X.~Hu, W.~Pan, A spatially adaptive high-order meshless method
  for fluid--structure interactions.
\newblock Computer Methods in Applied Mechanics and Engineering \textbf{355},
  67--93 (2019)

\bibitem{tolstykh2003using}
A.~Tolstykh, D.~Shirobokov, On using radial basis functions in a “finite
  difference mode” with applications to elasticity problems.
\newblock Computational Mechanics \textbf{33}(1), 68--79 (2003)

\bibitem{oanh2022approach}
D.T. Oanh, N.M. Tuong, An approach to adaptive refinement for the rbf-fd method
  for 2d elliptic equations.
\newblock Applied Numerical Mathematics \textbf{178}, 123--154 (2022)

\bibitem{toth2022h}
B.~T{\'o}th, A.~D{\"u}ster, h-adaptive radial basis function finite difference
  method for linear elasticity problems.
\newblock Computational Mechanics pp. 1--20 (2022)

\bibitem{fan2019adaptive}
L.~Fan, Adaptive meshless point collocation methods: investigation and
  application to geometrically non-linear solid mechanics.
\newblock Ph.D. thesis, Durham University (2019)

\bibitem{mishra2020adaptive}
P.K. Mishra, L.~Ling, X.~Liu, M.K. Sen, Adaptive radial basis function
  generated finite-difference (rbf-fd) on non-uniform nodes using $ p
  $-refinement.
\newblock arXiv preprint arXiv:2004.06319  (2020)

\bibitem{milewski2021higher}
S.~Milewski, Higher order schemes introduced to the meshless fdm in elliptic
  problems.
\newblock Engineering Analysis with Boundary Elements \textbf{131}, 100--117
  (2021)

\bibitem{albuquerque2021generalized}
A.~Albuquerque-Ferreira, M.~Ure{\~n}a, H.~Ramos, The generalized finite
  difference method with third-and fourth-order approximations and treatment of
  ill-conditioned stars.
\newblock Engineering Analysis with Boundary Elements \textbf{127}, 29--39
  (2021)

\bibitem{liszka1996hp}
T.~Liszka, C.~Duarte, W.~Tworzydlo, hp-meshless cloud method.
\newblock Computer Methods in Applied Mechanics and Engineering
  \textbf{139}(1-4), 263--288 (1996)

\bibitem{jancic_p_refined}
M.~Jančič, J.~Slak, G.~Kosec, in \emph{2021 6th International Conference on
  Smart and Sustainable Technologies (SpliTech)} (2021), pp. 01--06.
\newblock \doi{10.23919/SpliTech52315.2021.9566401}

\bibitem{duarte1996hp}
C.A. Duarte, J.T. Oden, An hp adaptive method using clouds.
\newblock Computer methods in applied mechanics and engineering
  \textbf{139}(1-4), 237--262 (1996)

\bibitem{slak2019generation}
J.~Slak, G.~Kosec, On generation of node distributions for meshless pde
  discretizations.
\newblock SIAM journal on scientific computing \textbf{41}(5), A3202--A3229
  (2019)

\bibitem{slak2021medusa}
J.~Slak, G.~Kosec, Medusa: A c++ library for solving pdes using strong form
  mesh-free methods.
\newblock ACM Transactions on Mathematical Software (TOMS) \textbf{47}(3),
  1--25 (2021)

\bibitem{depolli_parallel_2022}
M.~Depolli, J.~Slak, G.~Kosec, Parallel domain discretization algorithm for
  {RBF}-{FD} and other meshless numerical methods for solving {PDEs}.
\newblock Computers \& Structures \textbf{264}, 106,773 (2022).
\newblock Publisher: Elsevier

\bibitem{duh2021fast}
U.~Duh, G.~Kosec, J.~Slak, Fast variable density node generation on parametric
  surfaces with application to mesh-free methods.
\newblock SIAM Journal on Scientific Computing \textbf{43}(2), A980--A1000
  (2021)

\bibitem{wendland2004scattered}
H.~Wendland, \emph{Scattered Data Approximation}.
\newblock Cambridge Monographs on Applied and Computational Mathematics
  (Cambridge University Press, 2004).
\newblock \doi{10.1017/CBO9780511617539}

\bibitem{DAVYDOV2023115031}
O.~Davydov, D.T. Oanh, N.M. Tuong, Improved stencil selection for meshless
  finite difference methods in 3d.
\newblock Journal of Computational and Applied Mathematics \textbf{425},
  115,031 (2023).
\newblock \doi{https://doi.org/10.1016/j.cam.2022.115031}.
\newblock
  \urlprefix\url{https://www.sciencedirect.com/science/article/pii/S037704272200629X}

\bibitem{imex}
M.~Jančič, F.~Strniša, G.~Kosec, in \emph{2022 7th International Conference
  on Smart and Sustainable Technologies (SpliTech)} (2022), pp. 01--04.
\newblock \doi{10.23919/SpliTech55088.2022.9854342}

\bibitem{heuer2001hp}
N.~Heuer, M.E. Mellado, E.P. Stephan, hp-adaptive two-level methods for
  boundary integral equations on curves.
\newblock Computing \textbf{67}(4), 305--334 (2001)

\bibitem{bayona2019insight}
V.~Bayona, An insight into rbf-fd approximations augmented with polynomials.
\newblock Computers \& Mathematics with Applications \textbf{77}(9), 2337--2353
  (2019)

\bibitem{tominec}
I.~Tominec, E.~Larsson, A.~Heryudono, A least squares radial basis function
  finite difference method with improved stability properties.
\newblock SIAM Journal on Scientific Computing \textbf{43}(2), A1441--A1471
  (2021).
\newblock \doi{10.1137/20M1320079}.
\newblock \urlprefix\url{https://doi.org/10.1137/20M1320079}.
\newblock
  {\href{https://arxiv.org/abs/https://doi.org/10.1137/20M1320079}{{https://doi.org/10.1137/20M1320079}}}

\bibitem{bayona2019comparison}
V.~Bayona, Comparison of moving least squares and rbf+ poly for interpolation
  and derivative approximation.
\newblock Journal of Scientific Computing \textbf{81}(1), 486--512 (2019)

\bibitem{daniel2018adaptive}
P.~Daniel, A.~Ern, I.~Smears, M.~Vohral{\'\i}k, An adaptive hp-refinement
  strategy with computable guaranteed bound on the error reduction factor.
\newblock Computers \& Mathematics with Applications \textbf{76}(5), 967--983
  (2018)

\bibitem{eigenweb}
G.~Guennebaud, B.~Jacob, et~al.
\newblock Eigen v3.
\newblock \url{http://eigen.tuxfamily.org} (2010)

\bibitem{wakeni2021p}
M.F. Wakeni, A.~Aggarwal, L.~Kaczmarczyk, A.T. McBride, I.~Athanasiadis, C.J.
  Pearce, P.~Steinmann, A p-adaptive, implicit-explicit mixed finite element
  method for diffusion-reaction problems.
\newblock International Journal for Numerical Methods in Engineering  (2022)

\bibitem{kosec_weak_2019}
G.~Kosec, J.~Slak, M.~Depolli, R.~Trobec, K.~Pereira, S.~Tomar, T.~Jacquemin,
  S.P. Bordas, M.A. Wahab, Weak and strong from meshless methods for linear
  elastic problem under fretting contact conditions.
\newblock Tribology International \textbf{138}, 392--402 (2019).
\newblock Publisher: Elsevier

\bibitem{Slaughter_2002}
W.S. Slaughter, \emph{Three-Dimensional Problems} (Birkh{\"a}user Boston,
  Boston, MA, 2002), pp. 331--386.
\newblock \doi{10.1007/978-1-4612-0093-2_9}.
\newblock \urlprefix\url{https://doi.org/10.1007/978-1-4612-0093-2_9}

\end{thebibliography}


\end{document}